\documentclass[11pt,reqno]{amsart}

\usepackage[english]{babel} % Palabras por defecto en Inglés, como fechas, títulos 

\usepackage[pagewise]{lineno}%\linenumbers

\usepackage{anysize}
\usepackage{verbatim}
\marginsize{2cm}{2cm}{3cm}{3cm} % Márgenes izquierda, derecha,arriba,abajo

\usepackage{color}%escribir en colores 

\usepackage{tikz}%fugres
\usepackage{dsfont}
\usetikzlibrary{chains,arrows.meta,quotes,bending}
\usetikzlibrary{automata, positioning}
\usepackage{capt-of}%Borra el caption

\usepackage{mathrsfs}
\usepackage{subfig}
\usepackage{amssymb}

\usepackage{parskip}
\usepackage{enumerate}

\usepackage[colorlinks,citecolor=red,urlcolor=blue]{hyperref} %Pequeño cuadro para la referencia
\usepackage[numbers,sort]{natbib} %COMANDOS BIBLIOGRAFIA
\numberwithin{equation}{section}%numerar las ecuaciones de acuerdo a seccion.

\newtheorem{theorem}{Theorem}[section]
\newtheorem{lemma}[theorem]{Lemma}
\newtheorem{proposition}[theorem]{Proposition}

\theoremstyle{definition}
\newtheorem{definition}[theorem]{Definition}
\theoremstyle{remark}

\newtheorem{remark}[theorem]{\bf Remark}

\newcommand{\NN}{\mathbb{N}}

\newcommand{\RR}{\mathbb{R}}
\newcommand{\R}{\mathbb{R}}

\DeclareMathAlphabet{\mathdutchcal}{U}{dutchcal}{m}{n}

\newcommand{\CaixaPreta}{\vrule Depth0pt height5pt width5pt}
\newcommand{\bdemo}{\noindent {\bf \em Proof.} \hspace{2mm}}
\newcommand{\fdemo}{\hfill \CaixaPreta \vspace{3mm}}%COMANDOS DEMOSTRACION 

\begin{document}

\title[]{Asymptotic behavior of the fundamental solution of space-time fractional equations with a reaction term.}

\author{Luciano Abadias}
\address{Departamento de Matemáticas,
        Instituto Universitario de Matemáticas y Aplicaciones, Universidad de Zaragoza, Spain.}
\email{labadias@unizar.es}
\thanks{First author has been  partially supported by Grant PID2022-137294NBI00 funded by MICIU/AEI/10.13039/501100011033 and ``ERDF/EU'' and Project E48-20R, Gobierno de Aragón, Spain.} 

%%%%%%%%%%%%%%%%%%%%%%%%

\author[C. Carrasco]{Claudio Carrasco}
\address{Departamento de Matem\'aticas, Facultad de Ciencias, Universidad de Chile.}
\email{claudio.carrasco.g@ug.uchile.cl}
\thanks{Second author has been supported by Chilean research ANID-CONICYT for national doctoral 21240764, and partially supported by Chilean research grant FONDECYT 1221271. } 

%%%%%%%%%%%%%%%%%%%%%%%

\author[J.C. Pozo]{Juan C. Pozo}
\address{Departamento de Matem\'aticas, Facultad de Ciencias, Universidad de Chile.}
\email{jpozo@uchile.cl}
\thanks{Third author has been partially supported by Chilean research grant FONDECYT 1221271.} 

%%%%%%%%%%%%%%%%%%%%%%%%

\date{\today}

\maketitle

\thispagestyle{empty}

%%%%%%%%%%%%%%%%%%%%%%%%%%%%%%%%%%%%%%%%%%

\begin{abstract} 
In this paper, we consider a space-time fractional partial differential equation with a reactive term. We describe the speed of invasion of its fundamental solution, extending recent results in this topic, which had been proved for the one dimensional spatial setting and some fractional parameters involved in the equation. The key tool to achieve these results in a wide range of cases (any spatial dimension, any time and spatial fractional differential parameters, and any polynomial speed of invasion) is the subordination. %\textcolor{red}{Hay que tener en cuenta que las ecuaciones que estamos tomando son no locales tanto en espacio como en tiempo, así que yo cambiaría el título.}
\end{abstract}

\smallskip

%%%%%%%%%%%%%%%%%%%%%%%%%%%%%%%%%%%%%%%%%%
%%%%%%%%%%%%%%%%%%%%%%%%%%%%%%%%%%%%%%%%%%
%%%%%%%%%%%%%%%%%%%%%%%%%%%%%%%%%%%%%%%%%%

\section{Introduction and main results}

The main aim of this paper is the study of the speed of invasion of the fundamental solution of the following {\em non-local partial differential equation with a reactive term}:
\begin{equation}\label{Main:Equation}
\partial_t^\alpha u(t,x)+(-\Delta)^{\varrho} u(t,x)= u(t,x),\quad t>0,\ x\in\RR^d,
\end{equation}
where the operator $(-\Delta)^{\varrho}$ stands for the fractional laplacian of order $\varrho>0$, and $\partial_t^\alpha$ denotes the time-fractional derivative of order $\alpha\in(0,1)$ in thesense of Caputo, which for a smooth enough function, is defined by 
\[
\partial_t^\alpha f(t)=\frac{1}{\Gamma(1-\alpha)}\int_0^t \frac{f'(s)}{(t-s)^\alpha}\,ds.
\]

The left-side part of the equation \eqref{Main:Equation} represents the classical fully nonlocal diffusion partial differential equation. The time fractional derivative provides the subdiffusive behaviour of the solution, and the spatial fractional laplacian the superdiffusive one. In essence, a nonlocal diffusive partial differential equation is an evolution equation that combines nonlocal di\-ffu\-sion ope\-ra\-tors in the space variable with time-fractional derivatives (and generalizations). In recent years, the investigation into this type of partial differential equations has received considerable attention. This interest is primarily driven by their profound connections with non-local transport phenomena, the control of stochastic jump processes, the depiction of anomalous diffusion in physics, and the analysis of memory effects in parabolic equations. While not exhaustive, interested readers are referred to \cite{Ascione-Leonenko-Pirozzi-2020, Ascione-Mishura-Pirozzi-2022, Dier-Kemppainen-Siljander-Zacher-2020, Kemppainen-Siljander-Zacher-2017, Shlesinger-2017, Uchaikin-2000, Schneider-Wyss-1989, Kochubei-1989, Eidelman-Kochubei-2004, Cortazar-Quiros-Wolanski-2021:1, Cortazar-Quiros-Wolanski-2021:2}, and the references therein.

%%%%%%%%%%%%%%%%%%%%%%%%%%%%%%%%%%%%%%%

The prototype of subdiffusive partial differential equations is the so-called {\em time-fractional heat equation}
\begin{equation}\label{Heat:Equation:Memory}
\partial_t^\alpha u-\Delta u=0.
\end{equation}
For the well-posedness and regularity of the solutions of \eqref{Heat:Equation:Memory} in $\RR^d$, see, for example, the work of  Eidelman and Kochubei in \cite{Eidelman-Kochubei-2004}; Allen, Caffarelli and Vasseur in \cite{Allen-Caffarelli-Vasseur-2017}; and Kemppainen, Siljander and Zacher in \cite{Kemppainen-Siljander-Zacher-2017}. Optimal asymptotic decay rates of the solutions have been obtained by Kemppainen, Siljander, Vergara, and Zacher in \cite{Kemppainen-Siljander-Vergara-Zacher-2016}, and more recently by Cortazar, Quiroz, and Wolanski in \cite{Cortazar-Quiros-Wolanski-2021:1}. For the Dirichlet problem on bounded domains, Vergara and Zacher \cite{Vergara-Zacher-2017} show large time decay estimates of the solutions. Recently, Chan, Gomez-Castro, and V\'azquez \cite{Chan-Gomez-Vazquez-2024} studied the global well-posedness of \eqref{Heat:Equation:Memory} for singular solutions on bounded domains.

%%%%%%%%%%%%%%%%%%%%%%%%%%%%%%%%%%%%%

The previous paragraph makes clear that the study of this class of equations has reached a certain level of maturity. However, understanding the properties of solutions to subdiffusive partial differential equations, such as \eqref{Heat:Equation:Memory}, with reactive terms remains comparatively less developed, with many basic questions still unresolved. The inclusion of reactive terms introduces considerable technical challenges, creating a {\em competition} among the effects of the diffusion operator, the reactive terms, and the {\em memory} induced by the time-fractional derivative.  %In such a context, we highlight the works made by Dipierro, Pellacci, Valdinoci, and Verzini in \cite{Dipierro-Pellacci-Valdinoci-Verzini-2021} who studied the {\em speed of invasion} of the fundamental solution of the equation \eqref{Main:Equation}. 

%%%%%%%%%%%%%%%%%%%%%%%%%%%%%%%%%%%%

In such a context, we highlight the work made by Dipierro, Pellacci, Valdinoci, and Verzini \cite{Dipierro-Pellacci-Valdinoci-Verzini-2021} who have studied the so-called {\em speed of invasion} of the solution of the following equation
\begin{equation}
\partial_t^\alpha u(t,x)-\Delta u(t,x)=u(t,x),\quad t>0,\ x\in\RR^d,\qquad\qquad u(0,x)=\delta_0(x),\quad x\in\RR^d,\label{Main:Equation:2}
\end{equation}
for the one dimensional case $d=1,$ where $\delta_0$ denotes the Dirac delta distribution. Roughly speaking, the invasion speed of the fundamental solution of equation \eqref{Main:Equation:2} serves to differentiate between whether the solution exhibits convergent or divergent behavior as the spatial variable $x$ varies across space on spheres whose radius grows with time according to a function $\theta(t)$. More precisely, if  $u_{\alpha,1}(t,x)$ denotes the solution to equation \eqref{Main:Equation:2}, it becomes an intriguing problem to characterize the behavior of $u_{\alpha,1}(t,x)$ on regions of the form ${|x|=\theta(t)}$, where $\theta$ is an increasing unbounded function, and to analyze the asymptotics of $u_{\alpha,1}(t,\theta(t)\vec{v})$ for large times. In \cite{Dipierro-Pellacci-Valdinoci-Verzini-2021}, the authors considered the function $\theta(t) = mt^\beta$ for $m > 0$ and $\beta > 0$, and derived the following result.
 
\begin{theorem}\label{Theo:DPVV} Let $d=1$ and $\alpha\in[\frac{1}{2},1)$. Let $u_{\alpha,1}$ be the solution of \eqref{Main:Equation:2}, $\beta\in(0,\frac12)$ and $m>0$. Then
\[
u_{\alpha,1}(t,mt^\beta)\ge m\alpha t^{\beta-2\alpha} E_{\alpha}(t^\alpha)\big(1-o(1)\big).
\]
In particular, we have that 
\[
u_{\alpha,1}(t,mt^\beta) \to \infty \quad \text{as}\quad t\to\infty.
\]
\end{theorem} 

%%%%%%%%%%%%%%%%%%%%%%%%%%%%%%%%%%%

Theorem \ref{Theo:DPVV} seems to reveal an important physical feature of the subdiffusion processes modeled by time-fractional operators: namely, the memory effect of the fractional derivative \emph{slows down} the invasion rate compared to the classical case, in which diffusion occurs with a linear speed. This result establishes that such invasion still occurs, with a power-law time dependence as close as we wish to the square root function. Additionally, it opens up several very interesting research directions possed in \cite{Dipierro-Pellacci-Valdinoci-Verzini-2021}, which we address in this paper.

%%%%%%%%%%%%%%%%%%%%%%%%%%%%%%%%%%%

\begin{enumerate}[{\bf\em(Q1)}]
\item To find the {\em optimal speed} distinguishing the diverging from the vanishing behavior of the fundamental solution of \eqref{Main:Equation:2}, and determine the asymptotic behavior in the case $\beta\ge \frac12$.

\item To obtain a result in the spirit of Theorem \eqref{Main:Equation:2} for the highly nonlocal regime $\alpha\in(0,\frac12)$.

\item To understand the asymptotic behavior of the fundamental solution of \eqref{Main:Equation:2} in every spatial dimension $d\ge 2$, and not only the one-dimensional case.
\end{enumerate}

%%%%%%%%%%%%%%%%%%%%%%%%%%%%%%%%%%%

The primary objective of this paper is to answer the questions mentioned above and enhance the understanding of the speed of invasion of the fundamental solution of the equation \eqref{Main:Equation}. In such a context, we present one of our main results.

\begin{theorem}\label{Theorem:rho=2}
Let $d\in\NN$, $\alpha\in(0,1)$, $\gamma_\alpha:=(1-\alpha)\alpha^{\frac{\alpha}{1-\alpha}}$, $\displaystyle M_\alpha=\inf \left\{n\in\NN\ : \ n>\left(\frac{2}{\gamma_\alpha}\right)^{\frac{1-\alpha}{\alpha}}\right\}$ and $\vec{v}\in\mathbb{S}^{d-1}$. Let $u_{\alpha,1}$ be the solution of \eqref{Main:Equation:2}. Then, the following statements hold.
\begin{enumerate}[$(i)$]
\item If $\beta\in(0,1)$ and $m>0$ then 
\[
\lim_{t\to\infty}u_{\alpha,1}(t,mt^\beta \vec{v})=\infty.
\]

\item If $\beta>1$ and $m>0$ then %=\begin{cases}\infty, 
\[
\lim_{t\to\infty}u_{\alpha,1}(t,mt^\beta \vec{v})=0.
\]
\item If $\beta=1$ and $0<m<2\sqrt{1-\gamma_\alpha}$ then 
\[
\lim_{t\to\infty}u_{\alpha,1}(t,mt \vec{v})=\infty.
\]
\item If $\beta=1$ and $m>2M_\alpha\sqrt{1-\frac{\gamma_\alpha}{M_\alpha}}$, then 
\[
\lim_{t\to\infty}u_{\alpha,1}(t,mt \vec{v})=0.
\]
\end{enumerate}
\end{theorem}
\begin{remark} The proof of Theorem \ref{Theorem:rho=2} relies on two key components: the so-called \emph{subordination principle} in the sense of Pr\"uss (refer to, e.g., \cite[Chapter 4]{Pruss-1993} or \cite[Chapter 3]{Bajlekova-2001}), and the asymptotic expansions of the so-called \emph{Wright functions} (see, for instance, \cite[Appendix F]{Mainardi-2010}). Using the subordination principle, we derive a representation of the fundamental solution of \eqref{Main:Equation:2} in terms of the Wright functions, and by a delicate analysis, we demonstrate that the asymptotic expansion of the Wright functions governs the asymptotic behavior of $u_{\alpha,1}$. This approach differs entirely from the one employed in \cite{Dipierro-Pellacci-Valdinoci-Verzini-2021}, enabling us to address the three questions mentioned above comprehensively, and generalizing the results to another evolution equations.
\end{remark}

On the other hand, Theorem \ref{Theorem:rho=2} establishes that the speed of invasion of $u_{\alpha,1}$ does not depend on the parameter $\alpha$, except in the boundary case $\beta=1$. This observation leads to the next question: \emph{What is the influence of the diffusion operator?} In such a context, we study the following equation:
\begin{equation}\label{Main:Equation:3}
\partial_t^\alpha u(t,x)+(-\Delta)^{\varrho} u(t,x)= u(t,x),\quad t>0,\ x\in\RR^d,\qquad \qquad u(0,x)=\delta_0(x),
\end{equation}
where the operator $(-\Delta)^{\varrho}$ stands for the fractional Laplacian of order $\varrho>0$.

\begin{theorem}\label{Theo:rho:general}
Let $d\in\NN$, $\alpha\in(0,1),$ $\varrho>0$ and $\vec{v}\in\mathbb{S}^{d-1}$. Let $u_{\alpha,\varrho}(t,x)$ the solution to \eqref{Main:Equation:3}. The following statements hold. 
\begin{enumerate}[$(i)$]
\item If $\varrho\in(0,1)$, $m>0$ and $\beta>0$ then
\[
\lim_{t\to\infty}u_{\alpha,\varrho}(t,mt^\beta \vec{v})=\infty.
\]
\item If $\varrho\in(1,\infty)$, $m>0$ and $\beta>1$ then
\[
\lim_{t\to\infty}u_{\alpha,\varrho}(t,mt^\beta \vec{v})=0.
\]
\end{enumerate}
\end{theorem}

The statement $(i)$ of Theorem \ref{Theo:rho:general} shows that a power-type function does not identify the optimal rate distinguishing the diverging and vanishing behavior of the fundamental solution of equation \eqref{Main:Equation:3} when $\varrho \in (0,1)$. In such a context, we present the following result.

\begin{theorem}\label{Theorem:rho:in:(0,2)}
Let $d\in\NN$, $\alpha\in(0,1)$, and $\vec{v}\in\mathbb{S}^{d-1}$. Let $u_{\alpha,\varrho}$ be the solution of \eqref{Main:Equation:3} with $\varrho\in(0,1)$. Define $\theta_{m,\beta}(t) = (e^{m t^{\beta}} - 1)$ for some $m > 0$ and $\beta > 0$. The following statements hold.
\begin{enumerate}[$(i)$]
\item If $\beta\in(0,1)$ and $m>0$ then 
\[
\lim_{t\to\infty}u_{\alpha,\varrho}(t,\theta_{m,\beta}(t)\vec{v})=\infty.
\]

\item If $\beta>1$ and $m>0$ then %=\begin{cases}\infty, 
\[
\lim_{t\to\infty}u_{\alpha,\varrho}(t,\theta_{m,\beta}(t)\vec{v})=0.
\]
\item If $\beta=1$, and $0<m\le \frac{1-\gamma_\alpha}{d+2\varrho}$, where $\gamma_\alpha=(1-\alpha)\alpha^{\frac{\alpha}{1-\alpha}}$, then 
\[
\lim_{t\to\infty}u_{\alpha,\varrho}(t,\theta_{m,\beta}(t)\vec{v})=\infty.
\]
\item If $\beta=1$ and $m>\frac{1}{d+2\varrho}$, then 
\[
\lim_{t\to\infty}u_{\alpha,\varrho}(t,\theta_{m,\beta}(t)\vec{v})=0.
\]
\end{enumerate}
\end{theorem}

\begin{remark}
The method developed to prove Theorem \ref{Theo:rho:general} does not provide information about the asymptotic behavior of $u_{\alpha,\varrho}$ in the case to $\varrho>1$ and $0<\beta<1$. The following result gives a partial answer to this problem in the one-dimensional case and $\beta\in\Big(0,\frac{1}{2\varrho}\Big)$. To prove such result, we need to get some new upper and lower estimates for the Mittag-Leffler function in terms of its derivative (see Lemmatas \ref{Lemma:upper:estimate} and \ref{Lemma:lower:estimate}), which are novel and inherently interesting by themselves, and generalize those obtained in \cite{Dipierro-Pellacci-Valdinoci-Verzini-2021}.
\end{remark}

\begin{theorem}\label{Theorem:one:dimension} Let $\alpha\in(0,1)$ and $\varrho\ge1$. Let $u_{\alpha,\varrho}(t,x)$ the solution to \eqref{Main:Equation:3}. If $d=1$, $m>0$, and $\beta\in\Big(0,\frac{1}{2\varrho}\Big)$ then
\[
\lim_{t\to\infty}u_{\alpha,\varrho}(t,mt^\beta \vec{v})=\infty.
\]

\end{theorem}

\subsection*{Organization of the paper} The paper is organized as follows. In Section \ref{Preliminaries}, we give some preliminary concepts and known results that we need for our work. In particular, we recall some features and properties of the so-called {\it Mittag-Leffler} and {\it Wright functions}, which are crucial in our analysis. In Section \ref{Proof:of:Main:Results}, we provide the proof of our main results when the speed of invasion is of polynomial type. In Section \ref{exponetial:invasion:speed} we show that the {\em optimal speed} distinguishing the diverging from the vanishing behavior of the fundamental solution of \eqref{Main:Equation:3} when $\varrho\in(0,1)$ is of exponential type.

%%%%%%%%%%%%%%%%%%%%%%%%%%%%%%%%%%%%%%%%%%%%%%%%%

\section{Preliminaries.}\label{Preliminaries}

In this section we summarize the most important preliminary results, concepts and definitions that we use throughout the paper. Along this paper, the Fourier transform of $v\in \mathcal{S}(\RR^d)$ is denoted by
\[
\widetilde{v}(\xi)=\int_{\RR^d} e^{-i\xi\cdot x}v(x)\,dx,\quad  \xi\in\R^d,
\]
extended as usual to $\mathcal{S}'(\RR^d),$ where $\mathcal{S}(\RR^d)$ denotes the Schwartz space.

The Laplace transform of a subexponential function $f\colon (0,\infty)\to \RR$ is defined by
\[
\widehat{f}(\lambda)=\int_{0}^\infty e^{-\lambda t}f(t)\,dt,
\]
whenever the last integral is convergent.

\subsection*{Asymptotic notation}In the context of asymptotic analysis of functions, we utilize the following classical notation. Let \( f, g \colon (0, \infty) \to \RR \) two real functions. We say that \( f = o(g) \) as \( x \to \infty \) if
\begin{align*}
    \lim_{x \to \infty} \frac{f(x)}{g(x)} = 0.
\end{align*}

We use the notation \( f = O(g) \) as \( x \to \infty \) to indicate the existence of constants \( M > 0 \) and \( x_0 \in \mathbb{R} \) such that
\begin{align*}
    |f(x)| \leq M g(x), \quad \text{for all } x > x_0.
\end{align*} Also, in this case we write $f \lesssim g.$ 

Additionally, we write \( f \sim g \) as \( x \to \infty \) if
\begin{align*}
    \lim_{x \to \infty} \frac{f(x)}{g(x)} = 1.
\end{align*}

\subsection*{Mittag-Leffler function} Let $\alpha,\beta>0$. The two parameter Mittag-Leffler function $E_{\alpha,\beta}$ is defined by
\begin{equation}\label{Mittag-Leffler}
E_{\alpha,\beta}(z)=\sum_{n=0}^\infty \frac{z^n}{\Gamma(\alpha n+\beta)},\quad z\in\RR.
\end{equation} 
In the case $\beta=1$, we will simply write $E_\alpha(z)$. Clearly, when $\alpha=1$ and $\beta=1$, this function coincides with the exponential function. The Mittag-Leffler function is essential for the treatment of time-fractional evolution equations. There is a rich theory of Mittag-Leffler functions that generalizes semigroup theory, where the standard reference is the monograph by Mainardi \cite{Mainardi-2010}, see also \cite{Kilbas-Srivastava-Trujillo-2006,Podlubny-1999}.

Next two lemmas can be found in \textcolor{red}{...}, and they will be applied along the paper.

\begin{lemma} (Derivative of Mittag-Leffler function)
For all $\alpha>0$ we have that 
\[
\frac{d}{dz} E_{\alpha}(z)=\frac{1}{\alpha}E_{\alpha,\alpha}(z), \quad z\in\RR.
\]
\end{lemma}

\begin{lemma} (Asymptotic behavior of Mittag-Leffler function)\label{asymp:mittag} Let $\alpha\in(0,1)$, then 
\[
E_\alpha(z)= -\frac{1}{\Gamma(1-\alpha)z}+O(z^{-2})\quad \text{as} \quad z\to-\infty,
\]

\[
E_\alpha(z)=\frac{1}{\alpha}\exp\Big(z^\frac1\alpha\Big)+O(z^{-1})\quad \text{as} \quad z\to\infty,
\]
and
\[E_{\alpha,\alpha}(z)=\frac{1}{\alpha}z^{\frac{1-\alpha}{\alpha}}e^{z^{\frac{1}{\alpha}}}+O(z^{2}), \quad \text{as} \quad z \to \infty \]
\end{lemma}

The following results have been obtained by Dipierro, Pellacci, Valdinoci and Verzini, see \cite{Dipierro-Pellacci-Valdinoci-Verzini-2021}.

\begin{lemma} Let $\alpha\in[\frac12,1)$ and $r\ge0$. Then
\[
E_\alpha(r)\le \alpha E_{\alpha}'(r)+1-\frac{1}{\Gamma(\alpha)}+\left(\frac{1}{\Gamma(1+\alpha)}-\frac{1}{\Gamma(2\alpha)}\right)r.
\]
\end{lemma}

\begin{lemma} Let $\alpha\in [\frac12,1)$ and $r>0$. Then
\[
E_\alpha(r)\ge \frac{\alpha}{r} E_{\alpha}'(r)-\frac{1}{\Gamma(\alpha)r}+1-\frac{1}{\Gamma(2\alpha)}.
\]
\end{lemma}

In what follows we extend above results for any $\alpha\in(0,1)$, not only for $\alpha\in[1/2,1)$. Such estimates are key points in the proof of Theorem \ref{Theorem:one:dimension}.

\begin{lemma}\label{Lemma:upper:estimate}
Let $n \in\{2,3,\cdots\}$ and $\alpha \in [\frac{1}{n},1]$. Then
\begin{equation}\label{upper:estimate}
    E_{\alpha}(r)\leq \alpha E'_{\alpha}(r)+\sum_{k=0}^{\left[\frac{3n-2}{2}\right]-1}\gamma_{k}r^k,\quad r\ge 0,
\end{equation}
where
\[
\gamma_k=\frac{1}{\Gamma(1+\alpha k)}-\frac{1}{\Gamma(\alpha+\alpha k)}.
\]
 \end{lemma}
 \bdemo Let $r\ge 0$. We note that 
\[
E_\alpha(r)=\sum_{k=0}^{[\frac{3n-2}{2}]-1}\frac{r^{k}}{\Gamma(1+\alpha k)}+\sum_{k=[\frac{3n-2}{2}]}^{\infty}\frac{r^{k}}{\Gamma(1+\alpha k)}
\]
Since $\alpha \geq \frac{1}{n}$ and $n\ge 2$, it follows that
\[
1+\alpha k\ge\frac{3}{2},\quad\text{for all}\quad k \geq \left[\frac{3n-2}{2}\right].
\]
Therefore,

\begin{align*}    
E_{\alpha}(r)    &\leq \sum_{k=0}^{[\frac{3n-2}{2}]-1}\frac{r^{k}}{\Gamma(1+\alpha k)}+\sum_{k=[\frac{3n-2}{2}]}^{\infty}\frac{r^{k}}{\Gamma(\alpha+\alpha k)}\\
& =\sum_{k=0}^{[\frac{3n-2}{2}]-1}\frac{r^{k}}{\Gamma(1+\alpha k)}-\sum_{k=0}^{[\frac{3n-2}{2}]-1}\frac{r^{k}}{\Gamma(\alpha+\alpha k)}+E_{\alpha,\alpha}(r)\\ 
&=\sum_{k=0}^{[\frac{3n-2}{2}]-1}\left(\frac{1}{\Gamma(1+\alpha k)}-\frac{1}{\Gamma(\alpha+\alpha k)}\right)r^k+\alpha E'_{\alpha}(r) \\
    &=\alpha E'_{\alpha}(r)+\sum_{k=0}^{[\frac{3n-2}{2}]-1}\gamma_{k}r^k,
\end{align*}
where 
\[
\gamma_k:=\frac{1}{\Gamma(1+\alpha k)}-\frac{1}{\Gamma(\alpha+\alpha k)}.
\]
This proves our claim.
 \fdemo
 
 \begin{lemma} \label{Lemma:lower:estimate} Let $n \in \{2,3,\cdots\}$ and $\alpha \in [\frac{1}{n},1]$. Then    
    \begin{align*} \label{3.5}
         E_{\alpha}(r)\geq \frac{\alpha}{r^{n-1}}E_{\alpha}'(r)-\frac{1}{r^{n-1}}\sum_{k=0}^{[n-1+\frac{1}{2\alpha}]}\lambda_{k}r^k+\frac{1}{r^{n-1}}\sum_{k=n-1}^{[n-1+\frac{1}{2\alpha}]}\beta_{k}r^k,
    \end{align*}
where 
\[
\lambda_k=\frac{1}{\Gamma(\alpha+k\alpha)},\quad\text{and}\quad \beta_k=\frac{1}{\Gamma(1+k\alpha-\alpha(n-1))}.
\] 
\end{lemma}

\bdemo Let $r>0$. Note that  
    \begin{align*}
        r^{n-1}E_{\alpha}(r)=\sum_{k=0}^{\infty}\frac{r^{k+(n-1)}}{\Gamma(1+k\alpha)}=\sum_{k=n-1}^{\infty}\frac{r^{k}}{\Gamma(1+k\alpha-\alpha(n-1))}.
    \end{align*}
On the other hand, if $k \geq (n-1)+\frac{1}{2\alpha}$ we have that
\begin{align*}
    \alpha+k\alpha\geq 1-\alpha(n-1)+k\alpha \geq 1-\alpha(n-1)+\left((n-1)+\frac{1}{2\alpha}\right)\alpha=\frac{3}{2},
\end{align*}
which in turn implies that
\[
\Gamma(\alpha+k\alpha) \geq \Gamma(1-\alpha(n-1)+k\alpha), \quad \text{for all}\quad k\ge  (n-1)+\frac{1}{2\alpha}.
\]
In consequence we have that
\begin{align*}
     r^{n-1}E_{\alpha}(r)&=\sum_{k=n-1}^{\infty}\frac{r^{k}}{\Gamma(1+k\alpha-\alpha(n-1))}\\
    &=\sum_{k=n-1}^{[n-1+\frac{1}{2\alpha}]}\frac{r^{k}}{\Gamma(1+k\alpha-\alpha(n-1))}+\sum_{k=[n-1+\frac{1}{2\alpha}]+1}^{\infty}\frac{r^{k}}{\Gamma(1+k\alpha-\alpha(n-1))}\\
    &\geq \sum_{k=n-1}^{[n-1+\frac{1}{2\alpha}]}\frac{r^{k}}{\Gamma(1+k\alpha-\alpha(n-1))}+\sum_{k=[n-1+\frac{1}{2\alpha}]+1}^{\infty}\frac{r^{k}}{\Gamma(\alpha+k\alpha)}\\
    &= E_{\alpha, \alpha}(r)-\sum_{k=0}^{[n-1+\frac{1}{2\alpha}]}\frac{r^k}{\Gamma(\alpha+k\alpha)}+\sum_{k=n-1}^{[n-1+\frac{1}{2\alpha}]}\frac{r^{k}}{\Gamma(1+k\alpha-\alpha(n-1))}
\end{align*}

It follows from  \eqref{upper:estimate} that
\begin{align*}
    r^{n-1}E_{\alpha}(r)&\ge\alpha E'_{ \alpha}(r)-\sum_{k=0}^{[n-1+\frac{1}{2\alpha}]}\frac{r^k}{\Gamma(\alpha+k\alpha)}+\sum_{k=n-1}^{[n-1+\frac{1}{2\alpha}]}\frac{r^{k}}{\Gamma(1+k\alpha-\alpha(n-1))}\\
    &=\alpha E'_{ \alpha}(r)-\sum_{k=0}^{[n-1+\frac{1}{2\alpha}]}\lambda_{k}r^k+\sum_{k=n-1}^{[n-1+\frac{1}{2\alpha}]}\beta_{k}r^k,
\end{align*}
where 
\[
\lambda_k:=\frac{1}{\Gamma(\alpha+k\alpha)}\quad\text{and} \quad \beta_k:=\frac{1}{\Gamma(1+k\alpha-\alpha(n-1))}.
\]
Our conclusion follows directly after dividing by $r^{n-1}$.
\fdemo

\subsection*{Wright functions}Let $\alpha>-1$ and $\beta\in\RR$. The Wright function $W_{\alpha,\beta}$ is defined by 
\begin{equation}\label{Wright}
W_{\alpha,\beta}(z)=\sum_{n=0}^\infty \frac{z^n}{n!\Gamma(\alpha n+\beta)},\quad z\in\RR.
\end{equation}
Now we provide a concise overview of various properties associated with these functions. In particular, we show that they are strongly related to Mittag-Leffler functions. For a comprehensive summary with abstract and applied results, we refer the reader to \cite{Abadias-Miana-2015}, \cite[Appendix F]{Gorenflo-Kilbas-Mainardi-Rogosin-2020} and \cite[Section 3.2]{jin_2022} and the references therein. We recall some of such results in the following lines, that we will use along the paper.

\begin{proposition}\label{Laplace:Wright} Let $\alpha\in(0,1)$ and $\beta>0$. Then
\[
E_{\alpha,\beta}(z)=\int_0^\infty e^{zs} W_{-\alpha,\beta-\alpha}(-s)\,ds.
\]
In particular, for $\alpha\in(0,1)$ we have that 
\[
E_{\alpha}(-z)=\int_0^\infty e^{-zs} W_{-\alpha,1-\alpha}(-s)\,ds.
\]
\end{proposition}

\begin{proposition}(Derivative of Wright functions)
Let $\alpha>-1$ and $\beta>0$
\[
\frac{d}{dz} W_{\alpha,\beta}(z)=W_{\alpha,\alpha+\beta}(z),\quad z>0.
\]
In particular, for $\alpha\in(0,1)$ we have that
\begin{equation}\label{Diff:Wa}
\frac{d}{dz} W_{-\alpha,1}(z)=W_{-\alpha,1-\alpha}(z),\quad z>0.
\end{equation} 
\end{proposition}
\begin{proposition}\label{Moments}
    Let $\alpha \in (0,1)$ and $\nu>-1$. Then  
    \begin{align*}
        \int_{0}^{\infty} W_{-\alpha,1-\alpha}(-r) r^{\nu}  \, dr=\frac{\Gamma(\nu+1)}{\Gamma(\nu \alpha +1)}.
    \end{align*}
\end{proposition}

\begin{proposition}\label{Asymp:Wa}(Asymptotic behavior of Wright functions) Let $\nu\in(0,1)$ and $\mu>0$ then 
\begin{align*}
    W_{-\nu,\mu}(z)=Y^{\frac{1}{2}-\mu}e^{-Y}\left(\sum_{m=0}^{M-1}A_m Y^{-m}+O(|Y|^{-M})  \right), \quad \text{as} \quad |z| \to \infty,
\end{align*}
 where $Y=(1-\nu)(-\nu^{\nu}z)^{\frac{1}{1-\nu}}$ and $A_m$ are some real numbers. 
In particular, for $\alpha \in (0,1)$ we have that

\begin{align*}
    W_{-\alpha,1}(-t^{1-\alpha})= (\gamma_\alpha t)^{-\frac{1}{2}} e^{-\gamma_\alpha t}\left(A_0+O((\gamma_\alpha t)^{-1})\right), \quad \text{as} \quad t \to \infty,
\end{align*}
where $\gamma_\alpha=(1-\alpha)\alpha^{\frac{\alpha}{1-\alpha}}$ and $\displaystyle A_0=\frac{1}{\sqrt{2\pi}(1-\alpha)^\alpha \alpha^{2\alpha-1}}.$
\end{proposition}

Above result was obtained by Wright (\cite{Wright-1940}) by means of the so-called {\em smoothing Stokes’ discontinuities method}, and independently by Mainardi and Tomirotti (\cite{Mainardi-Tomirotti-1995}) who used the so-called {\em saddle point method}. 

\begin{remark}\label{Remark:Asymptotic:Behavior}
The asymptotic behavior of the Wright functions implies that for any $\alpha\in(0,1)$ there are $C_\alpha>0$ and $t^\ast>0$ such that 
\[
W_{-\alpha,1}(-t^{1-\alpha}) \ge C_\alpha t^{-\frac{1}{2}}e^{-\gamma_\alpha t}, \quad t > t^\ast.
\]

Analogously, for any $\alpha\in(0,1)$ there are $C'_\alpha>0$ and $t^{\ast\ast}>0$ such that
\[
W_{-\alpha,1}(-2t^{1-\alpha}) \ge C'_\alpha t^{-\frac{1}{2}}e^{-\gamma_\alpha 2^{\frac{1}{1-\alpha}}t}, \quad t > t^{\ast\ast}.
\]

\end{remark}

\begin{proposition}\label{Bound:Wa}
    Let $\alpha \in (0,1)$. There exist $\kappa>0$ and $\sigma>0$ such that
\[
0<W_{-\alpha,1-\alpha}(-t) \leq \kappa e^{-\sigma t^{\frac{1}{1-\alpha}}}, \quad \text{for all}\quad t>0.
\]
In particular, we have that 
\[
W_{-\alpha,1-\alpha}(-t)<\kappa,\quad\text{for all}\quad t>0.
\]

\end{proposition}

For a complete summary of the properties of the Wright functions and some applications, we refer the reader to \cite[Chapter 3]{jin_2022} and \cite[Appendix F]{Mainardi-2010}. In particular, the proof of Proposition \ref{Bound:Wa} can be found in \cite[Theorem 3.18]{jin_2022}.
\section{Proof of main results}\label{Proof:of:Main:Results}
In this section we prove our main results. Initially, we find a representation of the solution of \eqref{Main:Equation:3} in terms of Wright functions.

\begin{lemma}\label{Lemma:Subordination} Let $\alpha\in(0,1)$, $\varrho >0$ and $d\in\NN$. The solution of the equation \eqref{Main:Equation:3} is given by

\begin{equation}\label{Sub:Fractional:Lapliacian}
u_{\alpha, \varrho}(t,x)=t^{-\alpha}\int_0^\infty u_{1,\varrho}(s,x)W_{-\alpha,1-\alpha}(-t^{-\alpha}s)\,ds, \quad t>0,\ x\in\RR^d,
\end{equation}
where $u_{1,\varrho}$ is the fundamental solution of 
\begin{equation}\label{Classical:Equation}
\partial_t u(t,x)+(-\Delta)^{\varrho} u(t,x)=u(t,x),\quad t>0,\ x\in\RR^d,
\end{equation}
and $W_{-\alpha,1-\alpha}$ is the Wright function given in \eqref{Wright}.
\end{lemma}

\bdemo  Taking the Fourier transform into the both sides of the equation \eqref{Main:Equation:3}, we have that 
\[
\partial_t^\alpha \widetilde{u}(t,\xi)+|\xi|^{2\varrho} \widetilde{u}(t,\xi)=\widetilde{u}(t,\xi),\quad t>0,\ \xi\in\RR^d,\qquad\qquad \widetilde{u}(0,\xi)=1,\quad \xi\in\RR^d.
\]
Therefore, the Fourier transform of $u_{\alpha,\varrho}$ is given by
\[
\widetilde{u}_{\alpha,\varrho}(t,\xi)=E_{\alpha}\big((1-|\xi|^{2\varrho})t^\alpha\big),\quad t> 0,\quad \xi\in\RR^d,
\]
where $E_\alpha$ is the Mittag-Leffler function defined in \eqref{Mittag-Leffler}. It follows from Proposition \ref{Laplace:Wright} that 
\[
E_{\alpha}\big((1-|\xi|^{2\varrho})t^\alpha\big)=\int_0^\infty e^{(1-|\xi|^{2\varrho})t^\alpha s} W_{-\alpha,1-\alpha}(-s)\,ds.
\]
By a simple change of variables, we have that
\begin{align*}
\widetilde{u}_{\alpha,\varrho}(t,\xi)&=t^{-\alpha}\int_0^\infty e^{(1-|\xi|^{2\varrho})s} W_{-\alpha,1-\alpha}(-t^{-\alpha}s)\,ds.
\end{align*}
On the other hand, the Fourier transform of the fundamental solution  $u_{1,\varrho}$ of the equation 
\eqref{Classical:Equation} is given by $\widetilde{u}_{1,\varrho}(t,\xi)=\exp\big((1-|\xi|^{2\varrho})t\big)$ for $t>0$ and $\xi\in\RR^d$. Therefore, we have that
\[
\widetilde{u}_{\alpha,\varrho}(t,\xi)=t^{-\alpha}\int_0^\infty \widetilde{u}_{1,\varrho}(s,\xi) W_{-\alpha,1-\alpha}(-t^{-\alpha}s)\,ds.
\]
Our claim follows from the uniqueness of Fourier transform. 
\fdemo

\begin{remark}
Since $u_{1,\varrho}$  depends radially on the spatial variable, it follows from representation \eqref{Sub:Fractional:Lapliacian} that $u_{\alpha,\varrho}$ also depends radially on the spatial variable. Therefore, we can consider an increasing function $\theta(t)$ and analyze the asymptotic behavior of $u_{\alpha,\varrho}(t,\theta(t)\vec{v})$.
\end{remark}

%%%%%%%%%%%%%%%%%%%%%%%%%%%%%%%%%%%%%%%%%%%%%%%%%

Now we are in position to prove one of our main results.

\subsection*{Proof of Theorem \ref{Theorem:rho=2}} It is a well-known fact that when $\varrho=1$, then $u_{1,1}$ the solution of the equation \eqref{Classical:Equation} is given by 
\[
u_{1,1}(t,x)= \frac{e^{-\frac{\|x\|^2}{4 t}}}{(4\pi t)^{d/2}}e^{t}, \quad t>0,\quad x\in\RR^d.
\]
Thus, according to Lemma \ref{Lemma:Subordination}, the fundamental solution of \eqref{Main:Equation} is given by 
\begin{equation}\label{expresionU}
u_{\alpha,1}(t,x)=t^{-\alpha}\int_0^\infty e^{s}\ \frac{e^{-\frac{\|x\|^2}{4 s}}}{(4\pi s)^{d/2}}W_{-\alpha,1-\alpha}(-t^{-\alpha}s)\,ds, \quad t> 0, \quad x\in\RR^d.
\end{equation}
Now, we divide the proof into two main cases. 

First, we face the cases $(i)$ and $(iii)$ getting lower bounds for the expression \eqref{expresionU}. Since all the functions involved in the integral of the right hand side of \eqref{expresionU} are positive, we have that 
\[
u_{\alpha,1}(t,x)\ge \int_t^{2t} e^{s}  \frac{e^{-\frac{\|x\|^2}{4s}}}{(4\pi s)^{d/2}}\  t^{-\alpha}W_{-\alpha,1-\alpha}(-t^{-\alpha}s)\,ds,\qquad t>0, \ x\in\RR^d.
\]
Since the mapping $s \mapsto e^se^{-\frac{\|x\|^2}{4s}}$ is increasing and $s \mapsto \frac{1}{(4\pi s)^{d/2}}$ is decreasing, it follows that
\[
u_{\alpha,1}(t,x)\ge e^{t}  \frac{e^{-\frac{\|x\|^2}{4 t}}}{(8\pi t)^{d/2}}\int_t^{2t}\ t^{-\alpha}W_{-\alpha,1-\alpha}(-t^{-\alpha}s)\,ds,\qquad t>0,\ x\in\RR^d.
\]
Further, by \eqref{Diff:Wa} one gets 
\begin{align*}
u_{\alpha,1}(t,x)&\ge e^{t}  \frac{e^{-\frac{\|x\|^2}{4 t}}}{(8\pi t)^{d/2}}\Big(W_{-\alpha,1}(-t^{1-\alpha})-W_{-\alpha,1}(-2t^{1-\alpha})\Big)
\\
&= e^{t}  \frac{e^{-\frac{\|x\|^2}{4t}}}{(8\pi t)^{d/2}}W_{-\alpha,1}(-t^{1-\alpha})\left(1-\frac{W_{-\alpha,1}(-2t^{1-\alpha})}{W_{-\alpha,1}(-t^{1-\alpha})}\right),\qquad t>0,\ x\in\RR^d.
\end{align*}
By considering $x=mt^\beta\vec{v}$ with $m>0$, $\beta>0$ and $\vec{v}\in\mathbb{S}^{d-1}$, we have that 
\[
u_{\alpha,1}(t,mt^\beta \vec{v})\ge e^{t}  \frac{e^{-\frac{m^2t^{2\beta}}{4t}}}{(8\pi t)^{d/2}}W_{-\alpha,1}(-t^{1-\alpha})\left(1-\frac{W_{-\alpha,1}(-2t^{1-\alpha})}{W_{-\alpha,1}(-t^{1-\alpha})}\right),\qquad t>0.
\]

It follows from Proposition \ref{Asymp:Wa} that
\[
\lim_{t\to\infty} \frac{W_{-\alpha,1}(-2t^{1-\alpha})}{W_{-\alpha,1}(-t^{1-\alpha})}=0.
\]
%with \textcolor{red}{*}
%\[
% W_{-\alpha,1}(-t^{1-\alpha}) \sim (Ct)^{-\frac{1}{2}}e^{-\gamma_\alpha t},\quad\text{as}\quad t\to\infty,
%\]

%where recall that $\gamma_\alpha=(1-\alpha)\alpha^{\frac{\alpha}{1-\alpha}}\in (0,1)$ provided that $\alpha\in(0,1)$. 
Therefore, by Remark \ref{Remark:Asymptotic:Behavior} we get 
\begin{equation}\label{Asymp:u:alpha}
u_{\alpha,1}(t,mt^\beta \vec{v})\gtrsim  t^{-\frac{d+1}{2}}\exp\Big(t-\frac{m^2 t^{2\beta-1}}{4}-\gamma_{\alpha}t\Big)\left(1-o(1)\right),\qquad t\to \infty.
\end{equation}

Particularly, under assumptions in the case $(i),$ we have $2\beta-1<1,$ then from \eqref{Asymp:u:alpha} follows that
\[
u_{\alpha,1}(t,mt^\beta \vec{v}) \to \infty,\quad \text{as}\quad t\to\infty,
\]
for any $m>0,$ proving the assertion $(i)$.

In the case $(iii),$ $\beta=1$, and by \eqref{Asymp:u:alpha} we have that 
\[
u_{\alpha,1}(t,mt \vec{v})\gtrsim t^{-\frac{d+1}{2}}\exp\Big(\big(1-\frac{m^2 }{4}-\gamma_\alpha\big)t\Big)\left(1-o(1)\right),\qquad t\to \infty.
\]
Therefore, if $m<2\sqrt{1-\gamma_\alpha}$ we have that $u_{\alpha,1}(t,mt\vec{v})\to\infty$ as $t\to \infty$, which proves the assertion $(iii)$.

The second main case in the proof treats the assertions $(ii)$ and $(iv).$ In both items, we need to find accurate upper bounds of $u_{\alpha,1}(t,x)$ using the expression \eqref{expresionU}. For this purpose, we fix $\beta\geq 1$ and define $\delta=\frac{\beta+1}{2}$. Now we decompose $u_{\alpha,1}$ as follows
\begin{align*}
u_{\alpha,1}(t,x)=\underbrace{t^{-\alpha}\int_0^{t^\delta}e^{s}\ \frac{e^{-\frac{\|x\|^2}{4 s}}}{(4\pi s)^{d/2}}W_{-\alpha,1-\alpha}(-t^{-\alpha}s)\,ds}_{I_1(t,x)} + \underbrace{t^{-\alpha}\int_{t^\delta}^\infty e^{s}\ \frac{e^{-\frac{\|x\|^2}{4 s}}}{(4\pi s)^{d/2}}W_{-\alpha,1-\alpha}(-t^{-\alpha}s)\,ds.}_{I_2(t,x)}
\end{align*}

We analyze both terms separately in each one cases $(ii)$ and $(iv)$. We begin estimating $I_1(t,x)$. By Proposition \ref{Bound:Wa} we have that $W_{-\alpha,1-\alpha}(-t^\alpha s)< \kappa$ for some $\kappa>0$. Hence, 
\[
I_1(t,x)\leq \kappa t^{-\alpha}\int_0^{t^\delta}e^{s}\ \frac{e^{-\frac{\|x\|^2}{4 s}}}{(4\pi s)^{d/2}}\,ds \leq \kappa\frac{t^{-\alpha}e^{t^\delta}}{\pi^{d/2}}\int_0^{t^\delta} \frac{e^{-\frac{\|x\|^2}{4 s}}}{(4s)^{d/2}}\,ds,\quad t>0,\ x\in\R^d.
\]
By making the change of variable $\sigma=\frac{\|x\|^2}{4s}$, we obtain that
\[
I_1(t,x) \leq \kappa \frac{t^{-\alpha}e^{t^\delta}\|x\|^{2-d}}{4\pi^{d/2}}\int_{\frac{\|x\|^2}{4t^\delta}}^\infty e^{-\sigma}\sigma^{\frac{d}{2}-2}\,d\sigma=\kappa\frac{t^{-\alpha}e^{t^\delta}\|x\|^{2-d}}{4\pi^{d/2}}\Gamma\left(\frac{d}{2}-1,\frac{\|x\|^2}{4t^\delta}\right),\quad t>0,\ x\in\R^d,
\]
where $\Gamma(\cdot,\cdot)$ denotes the incomplete Gamma function. Considering $x=mt^\beta\vec{v}$ with $m>0$, $\vec{v}\in\mathbb{S}^{d-1}$, we have that
\[
I_1(t,mt^\beta\vec{v})\leq \kappa \frac{t^{-\alpha}e^{t^\delta}m^{2-d}t^{(2-d)\beta}}{4\pi^{d/2}}\Gamma\left(\frac{d}{2}-1,\frac{m^2t^{2\beta-\delta}}{4}\right),\quad t>0.
\]
It is well known (see, e.g., \cite[Section 2]{Temme-1996}) that 
\[
\Gamma(s,\tau)\sim \tau^{s-1}e^{-\tau}, \quad \text{as}\quad \tau\to\infty.
\]
Hence, we have that
\[
\frac{t^{-\alpha}e^{t^\delta}m^{2-d}t^{(2-d)\beta}}{4\pi^{d/2}}\Gamma\left(\frac{d}{2}-1,\frac{m^2t^{2\beta-\delta}}{4}\right)\sim \frac{4t^{-\alpha-\beta+1-(\beta+1)d/4}}{(4\pi)^{d/2}m^2}\exp\left(t^{\delta}-\frac{m^2t^{2\beta-\delta}}{4}\right),\quad \text{as}\quad t\to\infty.
\]

Therefore, in the case $(ii)$ ($\beta>1$) we have $\delta<\beta,$ and it follows easily by previous estimate that $I_1(t,mt^\beta\vec{v})\to 0$ as $t\to\infty.$ If we are under assumptions of the case $(iv),$ that is $\beta=1$ ($\delta=1$) and $m>2M_\alpha\sqrt{1-\frac{\gamma_\alpha}{M_\alpha}}$, with $\displaystyle M_\alpha=\inf \left\{n\in\NN\ : \ n>\left(\frac{2}{\gamma_\alpha}\right)^{\frac{1-\alpha}{\alpha}}\right\},$ it follows that $M_{\alpha}\geq 2,$ and then $m>2M_\alpha\sqrt{1-\frac{\gamma_\alpha}{M_{\alpha}}}>2.$  Then we also conclude that $I_1(t,mt\vec{v})\to 0$ as $t\to\infty.$

Finally, we focus on get an upper vanishing estimate for $I_2(t,x).$ Observe that

\begin{align*}
I_2(t,x)&\le \frac{t^{-\alpha-\frac{d\delta}{2}}}{(4\pi)^{\frac{d}{2}}}\int_{t^\delta}^\infty e^{s} e^{-\frac{\|x\|^2}{4 s}}\ W_{-\alpha,1-\alpha}(-t^{-\alpha}s)\,ds= \frac{t^{-\alpha-\frac{d\delta}{2}}}{(4\pi)^{\frac{d}{2}}}\sum_{n=1}^\infty \int_{nt^\delta}^{(n+1)t^\delta} e^{s}e^{-\frac{\|x\|^2}{4 s}}\ W_{-\alpha,1-\alpha}(-t^{-\alpha}s)\,ds\\
&\le \frac{t^{-\alpha-\frac{d\delta}{2}}}{(4\pi)^{\frac{d}{2}}}\sum_{n=1}^\infty  e^{(n+1)t^\delta}e^{-\frac{\|x\|^2}{4 (n+1)t^{\delta}}}\int_{nt^\delta}^{(n+1)t^\delta}W_{-\alpha,1-\alpha}(-t^{-\alpha}s)\,ds\\
&\leq \frac{t^{-\frac{d\delta}{2}}}{(4\pi)^{\frac{d}{2}}} \sum_{n=1}^\infty  e^{(n+1)t^\delta}e^{-\frac{\|x\|^2}{4 (n+1)t^{\delta}}} W_{-\alpha,1}(-nt^{\delta-\alpha}),\quad t>0,\ x\in\R^d,
\end{align*}
where in the last inequality we have applied the relation \eqref{Diff:Wa}. Considering $x=mt^\beta\vec{v}$ with $m>0$, $\vec{v}\in\mathbb{S}^{d-1}$, we have that 

\begin{equation*}\label{Upper:bound}
I_2(t,mt^\beta\vec{v})\le \frac{t^{-\frac{d\delta}{2}}}{(4\pi)^{\frac{d}{2}}}\sum_{n=1}^\infty e^{(n+1)t^\delta}e^{-\frac{m^2 t^{2\beta-\delta}}{4 (n+1)}} W_{-\alpha,1}(-nt^{\delta-\alpha}),\quad t>0.
\end{equation*}
Applying again Proposition \ref{Asymp:Wa}, we have that
\[
 W_{-\alpha,1}(-y^{1-\alpha}) \sim (Cy)^{-\frac{1}{2}}e^{-\gamma_\alpha y},\quad\text{as}\quad y\to\infty,
\]
with $C>0,$ which implies that for all $\varepsilon>0,$ there is $y_0>0$ such that $$\left|\frac{W_{-\alpha,1}(-y^{1-\alpha})}{(Cy)^{-\frac{1}{2}}e^{-\gamma_\alpha y}}-1\right|<\varepsilon,\quad y>y_0.$$ In particular, if $y=n^{\frac{1}{1-\alpha}}t^{\frac{\delta-\alpha}{1-\alpha}}$, it follows that $$\frac{W_{-\alpha,1}(-nt^{\delta-\alpha})}{C^{-\frac{1}{2}}e^{-\gamma_{\alpha}n^{\frac{1}{1-\alpha}}t^{\frac{\delta-\alpha}{1-\alpha}}}}\leq
 \frac{W_{-\alpha,1}(-nt^{\delta-\alpha})}{(Cn^{\frac{1}{1-\alpha}}t^{\frac{\delta-\alpha}{1-\alpha}})^{-\frac{1}{2}}e^{-\gamma_\alpha n^{\frac{1}{1-\alpha}}t^{\frac{\delta-\alpha}{1-\alpha}}}}\leq (1+\varepsilon)$$ for all $t>t_0:=\max\{y_0^{\frac{1-\alpha}{\delta-\alpha}},1\}$ and $n\in\NN.$ Then, for all $t>t_0,$ one gets $$I_2(t,mt^\beta\vec{v})\lesssim \frac{\,t^{-\frac{d\delta}{2}}}{(4\pi)^{\frac{d}{2}}}\sum_{n=1}^\infty \exp\left(-n\Big(\gamma_\alpha n^{\frac{\alpha}{1-\alpha}}t^{\frac{\delta-\alpha}{1-\alpha}}+\frac{m^2t^{2\beta-\delta}}{4 n(n+1)}-\Big(1+\frac{1}{n}\Big)t^\delta\Big)\right).$$
At this point, we distinguish the cases $(ii)$ and $(iv).$

For the case $(ii)$ we have that $$\gamma_\alpha n^{\frac{\alpha}{1-\alpha}}t^{\frac{\delta-\alpha}{1-\alpha}}+\frac{m^2t^{2\beta-\delta}}{4 n(n+1)}-\Big(1+\frac{1}{n}\Big)t^\delta\geq \gamma_\alpha t^{\frac{\delta-\alpha}{1-\alpha}}-2t^\delta=:g(t),\quad t>t_0,\ n\in\NN.$$ Since $\delta>1,$ $g(t)\to\infty$ as $t\to\infty,$  so $$I_2(t,mt^\beta\vec{v})\lesssim \frac{\,t^{-\frac{d\delta}{2}}}{(4\pi)^{\frac{d}{2}}}\sum_{n=1}^\infty e^{-ng(t)}=\frac{\,t^{-\frac{d\delta}{2}}}{(4\pi)^{\frac{d}{2}}} \frac{1}{e^{g(t)}-1}\to 0,\ \text{ as }t\to\infty.$$

In the case $(iv),$ recall that $\beta=1$ and $m>2M_\alpha\sqrt{1-\frac{\gamma_\alpha}{M_\alpha}}>2,$ where 
\[
\displaystyle M_\alpha=\inf \left\{n\in\NN\ : \ n>\left(\frac{2}{\gamma_\alpha}\right)^{\frac{1-\alpha}{\alpha}}\right\}.
\]
In such a context, we have that
\begin{align*}
I_2(t,mt \vec{v})\lesssim&\frac{ t^{-\frac{d}{2}}}{(4\pi)^{\frac{d}{2}}}\sum_{n=1}^{M_\alpha-1} \exp\left(-nt\Big(\gamma_\alpha n^{\frac{\alpha}{1-\alpha}}+\frac{m^2}{4 n(n+1)}-\Big(1+\frac{1}{n}\Big)\Big)\right)\\
&+\frac{t^{-\frac{d}{2}}}{(4\pi)^{\frac{d}{2}}}\sum_{n=M_{\alpha}}^{\infty} \exp\left(-nt\Big(\gamma_\alpha n^{\frac{\alpha}{1-\alpha}}+\frac{m^2}{4 n(n+1)}-\Big(1+\frac{1}{n}\Big)\Big)\right),\quad t>t_0.
\end{align*}
On one hand, it follows easily that $\gamma_\alpha n^{\frac{\alpha}{1-\alpha}}+\frac{m^2}{4 n(n+1)}-\Big(1+\frac{1}{n}\Big)\geq 0$ for all $1\leq n\leq M_{\alpha}-1,$ therefore the first above term in the upper estimate of $I_2(t,mt \vec{v})$ vanishes as $t\to\infty.$ On the other hand, we have that
\[
\Big(\gamma_\alpha n^{\frac{\alpha}{1-\alpha}}+\frac{m^2 }{4n(n+1)}-\Big(1+\frac{1}{n}\Big)\Big)t\ge \Big(\gamma_\alpha M_\alpha^{\frac{\alpha}{1-\alpha}}-2\Big)t=:h(t),\quad t>t_0,\ n\geq M_{\alpha}.
\]
Since $\gamma_\alpha M_\alpha^{\frac{\alpha}{1-\alpha}}-2>0$ one gets 
\[
\frac{t^{-\frac{d}{2}}}{(4\pi)^{\frac{d}{2}}}\sum_{n=M_{\alpha}}^{\infty} \exp\left(-nt\Big(\gamma_\alpha n^{\frac{\alpha}{1-\alpha}}+\frac{m^2}{4 n(n+1)}-\Big(1+\frac{1}{n}\Big)\Big)\right)\le \frac{t^{-\frac{d}{2}}}{(4\pi)^{\frac{d}{2}}} \frac{e^{-h(t)M_\alpha}}{1-e^{-h(t)}}\to 0,\text{ as }t\to\infty,
\]
which proves our claim.
\fdemo

\begin{remark} We note that the speed of invasion of \( u_{\alpha,1} \) does not depend on the parameter \( \alpha \), except in the boundary case where \( \beta = 1 \). In such a case, we have that 

\[
\lim_{t \to \infty} u_{\alpha,1}(t, mt \vec{v}) = 
\begin{cases}
\infty, & m < 2\sqrt{1 - \gamma_\alpha}, \\
0, & m \geq 2M_\alpha \sqrt{1 - \frac{\gamma_\alpha}{M_\alpha}}.
\end{cases}
\]

While our method does not provide information about the asymptotic behavior of \( u_{\alpha,1}(t, mt \vec{v}) \) for \( m \) within the range 

\[
2\sqrt{1 - \gamma_\alpha} \leq m \leq 2M_\alpha \sqrt{1 - \frac{\gamma_\alpha}{M_\alpha}},
\] 
this gap closes as \( \alpha \to 1 \). To clarify this, we observe that 
\[
\lim_{\alpha \to 1^-} \gamma_\alpha = 0 \quad \text{and} \quad \lim_{\alpha \to 1^-} M_\alpha = 1,
\]
which in turn implies 

\[
\lim_{\alpha \to 1^-} 2\sqrt{1 - \gamma_\alpha} = 2 \quad \text{and} \quad \lim_{\alpha \to 1^-} 2M_\alpha \sqrt{1 - \frac{\gamma_\alpha}{M_\alpha}} = 2.
\]

Thus, in the limiting case, we recover the classical result:

\[
\lim_{t \to \infty} u_{1,1}(t, mt^\beta \vec{v}) = \lim_{t \to \infty} \frac{1}{(4\pi t)^{\frac{d}{2}}} e^{t - \frac{m^2}{4}t} = 
\begin{cases}
\infty, & m < 2, \\
0, & m \geq 2.
\end{cases}
\] 

This observation confirms that the classical behavior of the heat equation is indeed reflected in our analysis as \( \alpha \) approaches to one.
\end{remark}

% In such a case, there is a full range

%Althouht, our method does not provide information about the asymptotic behavior of $u_\alpha(t,mt\vec{v})$. Nevertheless, since \[
%\lim_{\alpha\to1^-}\gamma_\alpha=0 \quad \text{and}\quad \lim_{\alpha\to1^-}M_\alpha=1,
%\]
%we have that 

%This implies that the results of the classical case can be recovered from our results. \textcolor{red}{*} \textcolor{red}{No entiendo esta ultima frase} Indeed, if $\alpha=1$ and $|x|=mt$, then
%\begin{align*}
%    \lim_{t \to \infty}u_{1}(t,mt^\beta \vect{e}) &=\lim_{t \to \infty}\frac{1}{(4\pi t)^{\frac{n}{2}}}e^{ t-\frac{m^2}{4}t} \\
%    &=\begin{cases}\infty,& m <2\\ 0,&m \geq 2\end{cases}.
%\end{align*}
%\end{remark}

\subsection*{Proof of Theorem \ref{Theo:rho:general}} By Lemma \ref{Lemma:Subordination}, we know that    
\[
u_{\alpha,\varrho}(t,x)=t^{-\alpha}\int_{0}^{\infty}u_{1,\varrho}(s,x)W_{-\alpha,1-\alpha}(-st^{-\alpha})\, ds, \quad t>0, \, x \in \RR^n,
\]
where $u_{1,\varrho}(t,x)$ is the solution of \eqref{Main:Equation:3} . Assume first that $\varrho\in(0,1)$. It follows from \cite[Chapter 3]{Jacob-2001} that there is a convolution semigroup $(\mu^\varrho_t)_{t\ge 0}$ such that 
\[
u_{1,\varrho}(t,x)=e^{t}\mu^\varrho_t(x).
\]
It is well known that $\mu_t^\varrho(x)$ is non-negative. Therefore, we have that  
\[
u_{\alpha,\varrho}(t,x) \geq t^{-\alpha}\int_{t}^{2t}e^sW_{-\alpha,1-\alpha}(-st^{-\alpha})\mu^{\varrho}_s(x)\, ds.  
\]

Despite the absence of explicit formulas for $\mu_t^\varrho(x)$, the asymptotic behavior is well understood, see e.g., \cite[Section 2]{vazquez-2017}. More precisely, there are positive constants $C_1$ and $C_2$ such that 
\[
C_1\frac{t}{(\|x\|^2 +t^\frac{1}{\varrho})^{\frac{d+2\varrho}{2}}}\le \mu^\varrho_t(x)\le C_2\frac{t}{(\|x\|^2 +t^\frac{1}{\varrho})^{\frac{d+2\varrho}{2}}},\quad t>0, \ x\in\RR^d.
\]
This implies that 
\[
u_{\alpha,\varrho}(t,x) \geq C_1 t^{-\alpha}\int_{t}^{2t}e^sW_{-\alpha,1-\alpha}(-st^{-\alpha})\, \frac{s}{(\|x\|^2 +s^\frac{1}{\varrho})^{\frac{d+2\varrho}{2}}}\,ds.  
\]

Since the mapping $s \mapsto se^s$ is increasing and $s \mapsto \frac{1}{(\|x\|^2+s^{\frac{1}{\varrho}})^{\frac{d+2\varrho}{2}}}$ is decreasing, it follows that
\begin{align*}
    u_{\alpha,\varrho}(t,x) &\geq C_1 t^{-\alpha}\frac{ t\, e^t }{((2t)^{\frac{1}{\varrho}}+\|x\|^2)^{\frac{d+2\varrho}{2}}} \int_{t}^{2t}W_{-\alpha,1-\alpha}(-st^{-\alpha})\, ds\geq C_1 \frac{t\, e^t}{((2t)^{\frac{1}{\varrho}}+\|x\|^2)^{\frac{d+2\varrho}{2}}} W_{-\alpha,1}(-t^{1-\alpha}),
\end{align*}
where in the last inequality we have applied \eqref{Diff:Wa}.
By considering $x=mt^\beta\vec{v}$ with $m>0$, $\beta>0$ and $\vec{v}\in\mathbb{S}^{d-1}$ one gets
 
 \[
u_{\alpha,\varrho}(t,mt^\beta \vec{v})\ge C_1 \frac{t\,e^t}{((2t)^{\frac{1}{\varrho}}+m^2t^{2\beta})^{\frac{d+2\varrho}{2}}}W_{-\alpha,1}(-t^{1-\alpha})\qquad t>0.
\]
Since $\gamma_\alpha\in(0,1)$ and
\[
 W_{-\alpha,1}(-t^{1-\alpha}) \sim (Ct)^{-\frac{1}{2}}e^{-\gamma_\alpha t},\quad\text{as}\quad t\to\infty,
\]
with $C$ a positive constant, we have that $u_{\alpha,\varrho}(t,mt^\beta \vec{v})\to \infty$ as $t\to \infty$, which proves the assertion $(i)$.

In order to prove the assertion $(ii)$, assume that both $\varrho,\beta>1$ and we decompose $u_{\alpha,\varrho}$ as follows
\begin{align*}
u_{\alpha,\varrho}(t,x)=t^{-\alpha}\int_0^{t^\delta}u_{1,\varrho}(s,x)W_{-\alpha,1-\alpha}(-t^{-\alpha}s)\,ds+ t^{-\alpha}\int_{t^\delta}^\infty u_{1,\varrho}(s,x)W_{-\alpha,1-\alpha}(-t^{-\alpha}s)\,ds,
\end{align*}
where $\delta=\frac{\beta+1}{2}$. Therefore, by considering $x=mt^\beta\vec{v}\,$ with $m>0$ and $\vec{v}\in\mathbb{S}^{d-1}$, we have
\begin{align*}
    |u_{\alpha,\varrho}(t,mt^\beta \vec{v})| \leq I_1(t)+I_2(t),
\end{align*}
where
\[
I_1(t)=t^{-\alpha}\int_{0}^{t^{\delta}} W_{-\alpha,1-\alpha}(-st^{-\alpha})|u_{1,\varrho}(s,mt^\beta \vec{v})|\, ds, 
\]
and
\[
I_2(t)=t^{-\alpha}\int_{t^{\delta}}^{\infty} W_{-\alpha,1-\alpha}(-st^{-\alpha})|u_{1,\varrho}(s,mt^\beta \vec{v})|\, ds.
\]

We analyze both terms separately. We start estimating $I_1(t)$. It is a well known fact that 
\[
u_{1,\varrho}(t,x)=e^{t}\nu^\varrho_t(x),
\] 
where $\nu^{\varrho}_t$ is the unique function whose Fourier transform is given by $\widehat{\nu^{\varrho}_t}(\xi)=\exp(-t|\xi|^{2\varrho})$ for $\xi\in\RR^d$. It is well known (see, e.g., \cite[Section 2]{Galaktionov-Pohozaev-2002}) that 
\[
\nu^{\varrho}_t(x)=t^{-\frac{d}{2\varrho}} F(t^{-\frac{1}{2\varrho}}x),\quad t>0, \ x\in\RR^d,
\]
where $\displaystyle F(y)=(2\pi)^{-\frac{d}{2}}\int_0^\infty e^{-s^{2\varrho}}(s\|y\|)^{d/2}J_{(d-2)/2}(s\|y\|)ds$ for $y\in\R^d$  and $J_{(d-2)/2}$ denotes the Bessel functions of first kind.  %\textcolor{red}{*} \textcolor{red}{Yo pondría en vez de se puede probar, es conocido... y daria referencia. Ademas, creo que esta formula aparece en el Blumental yno es exactamente asi, habria que revisarla, aunque en el paper que aparece a continuacion, 12 si la pone asi}.

We point out that there is not closed formulas for the function $F$. However, it has been proved in \cite[Proposition 2.1]{Galaktionov-Pohozaev-2002} that for any dimension $d$ there are constants $K>1$ and $\omega>0$ (depending on $d$ and $\varrho$) such that 
\[
|F(y)|\le K e^{-\omega \|y\|^{\frac{2\varrho}{2\varrho-1}}},\quad y\in\RR^d,
\]
% \textcolor{red}{aqui imagino que se usa otra vez que $W_{-\alpha,1-\alpha}(-st^{-\alpha})\leq 1,$ habria que dar referencia wue no encuentro la primera vez que se usa o incluso ponerlo en las propiedades de la seccion 2, y aqui referenciarlo }
where, $K=D\int_{\RR^d}e^{-q(1+\|y\|^2)^{\frac{\varrho}{2\varrho-1}}}\, dy$, for a some constants $D>1$ and $q>0$. Hence, the previous considerations together Proposition \ref{Bound:Wa} implies that %\textcolor{red}{*} 
\begin{align*}
I_1(t)&\leq C t^{-\alpha}e^{t^\delta}\int_{0}^{t^{\delta}} s^{-\frac{d}{2\varrho}} |F(mt^\beta s^{-\frac{1} {2\varrho}}\vec{v})|\, ds, \\
&\leq C  t^{-\alpha}e^{t^\delta}\int_{0}^{t^{\delta}} s^{-\frac{d}{2\varrho}} e^{-\omega(mt^{\beta} s^{-\frac{1}{2\varrho}})^{\frac{2\varrho}{2\varrho-1}}}ds,\\
&\leq C  t^{-\alpha}e^{t^\delta} \int_{0}^{t^\delta}s^{-\frac{d}{2\varrho}}e^{-\gamma s^{-\frac{1}{2\varrho-1}}}\, ds, \quad t\to\infty,
\end{align*}
where $\gamma=\gamma(m,t,\varrho,\beta):=\omega (mt^{\beta})^{\frac{2\varrho}{2\varrho-1}}$ and $C$ denotes a positive constant different from line to line and independent on $t.$
By making the following change of variables $u=\gamma s^{-\frac{1}{2\varrho -1}}$, we note that 
\begin{align*}
\int_{0}^{t^\delta}s^{-\frac{d}{2\varrho}}e^{-\gamma s^{-\frac{1}{2\varrho-1}}}\, ds&=(2\varrho-1)\gamma^{(2\varrho-1)(1-\frac{d}{2\varrho})}\int_{\gamma t^{-\frac{\delta}{2\varrho -1}}}^{\infty} u^{d\left(\frac{2\varrho-1}{2\varrho}\right)-2\varrho}e^{-u}\, du \\
&=(2\varrho-1)\gamma^{(2\varrho-1)(1-\frac{d}{2\varrho})} \Gamma\left(d\left(\frac{2\varrho-1}{2\varrho}\right)-2\varrho+1,\gamma t^{-\frac{\delta}{2\varrho-1}}\right)
\end{align*}
In consequence we have that
\begin{align*}
I_1(t)\leq Ct^{-\alpha}e^{t^\delta}\gamma^{(2\varrho-1)(1-\frac{d}{2\varrho})} \Gamma\left(d\left(\frac{2\varrho-1}{2\varrho}\right)-2\varrho+1,\,\omega m^{\frac{2\varrho}{2\varrho-1}} t^{\frac{2\varrho\beta-\delta}{2\varrho-1}}\right),\quad t>0,
\end{align*}
with $C>0$ independent on $t.$ Note that $t^{\frac{2\varrho\beta-\delta}{2\varrho-1}}\to\infty$ as $t\to\infty,$ since $\beta>\delta$ and $\varrho>1,$ so by using the asymptotic behavior of the incomplete Gamma function (see for example \cite[Section 2]{Temme-1996}) we have that 
$$
\Gamma\left(d\left(\frac{2\varrho-1}{2\varrho}\right)-2\varrho+1,\,\omega m^{\frac{2\varrho}{2\varrho-1}} t^{\frac{2\varrho\beta-\delta}{2\varrho-1}}\right)
\sim (\omega m^{\frac{2\varrho}{2\varrho-1}} t^{\frac{2\varrho\beta-\delta}{2\varrho-1}})^{\tau-1}e^{-\omega m^{\frac{2\varrho}{2\varrho-1}}t^{\frac{2\varrho \beta - \delta}{2\varrho-1}}},\quad t\to\infty,$$
%&\sim C (2\varrho-1)m^\varrho t^{-\alpha-\frac{d}{2\varrho}+2\varrho\beta}\, \ell(t)^{\delta} \ e^{-\ell(t)}e^{t^\delta}, \quad \text{as}\quad t\to\infty.
where $\tau:=d\left(\frac{2\varrho-1}{2\varrho}\right)-2\varrho+1$. Hence,
\begin{align*}
    I_{1}(t) \lesssim t^{-\alpha} (\omega m^{\frac{2\varrho}{2\varrho-1}} t^{\frac{2\varrho\beta-\delta}{2\varrho-1}})^{\tau-1}e^{-\omega m^{\frac{2\varrho}{2\varrho-1}}t^{\frac{2\varrho \beta - \delta}{2\varrho-1}}+t^{\delta}},\quad t\to\infty.
\end{align*}
Since $\delta<\beta$, it follows that 
\[
t^\delta-\omega m^{\frac{2\varrho}{2\varrho-1}}t^{\frac{2\varrho \beta - \delta}{2\varrho-1}}\to -\infty\quad\text{as}\quad t\to\infty,
\]
which implies that $I_1(t)\to 0$ as $t\to\infty$.

On the other hand,
\begin{align*}
    I_2(t)&=t^{-\alpha}\int_{t^{\delta}}^{\infty} W_{-\alpha,1-\alpha}(-st^{-\alpha}) e^s |\nu^\varrho_s(mt^\beta\vec{v})|\, ds \\
    & \leq C t^{-\alpha}\int_{t^{\delta}}^{\infty} e^s e^{-\tilde{m} s^{-\frac{1}{2\varrho -1}}t^{\frac{2\beta\varrho}{2\varrho-1}}}W_{-\alpha,1-\alpha}(-st^{-\alpha})\, ds
\end{align*}
for all $t>1,$ with $C$ a positive constant independent on $t$ and  $\Tilde{m}:=\omega m^{\frac{2\varrho}{2\varrho -1}}$.  Proceeding analogously to the proof of Theorem \ref{Theorem:rho=2} (ii), it follows that
\begin{align*}
    I_2(t) \lesssim  \sum_{n=1}^{\infty}\exp \left( -n\Big(\gamma_{\alpha}n^{\frac{\alpha}{1-\alpha}}t^{\frac{\delta-\alpha}{1-\alpha}}-\left(1+\frac{1}{n}\right) t^\delta+\frac{\Tilde{m}t^{\frac{2\varrho \beta-\delta}{2\varrho -1}}}{n(n+1)^{\frac{1}{2\varrho -1}}}\Big) \right ).
\end{align*}
 
Since $\beta>1$, we have that $\frac{\delta-\alpha}{1-\alpha}>\delta,$ so there is $t_0>1$ such that
\begin{align*}
    \gamma_{\alpha}n^{\frac{\alpha}{1-\alpha}}t^{\frac{\delta-\alpha}{1-\alpha}}-\left(1+\frac{1}{n}\right) t^\delta+\frac{\Tilde{m}t^{\frac{2\varrho \beta-\delta}{2\varrho -1}}}{n(n+1)^{\frac{1}{2\varrho -1}}}> \gamma_{\alpha}n^{\frac{\alpha}{1-\alpha}}t^{\frac{\delta-\alpha}{1-\alpha}}-\left(1+\frac{1}{n}\right) t^\delta>0
\end{align*}
for all $n\in\mathbb{N}$.  In the same way as in Theorem \ref{Theorem:rho=2} (ii)
\begin{align*}
    I_2(t) \to 0, \quad \text{as} \quad t \to \infty.
\end{align*}
\fdemo

\subsection*{Proof of Theorem \ref{Theorem:one:dimension}} We point out that when $\varrho>1$ the fundamental solution $u_{1,\varrho}$ of 
\[
\partial_t u(t,x)+(-\Delta)^{\varrho} u(t,x)=u(t,x),
\]
is sign changing. This is a big obstacle to apply in a sharp way the techniques developed above to study the speed of invasion of \eqref{Main:Equation:3}, for example in the case of $u_{\alpha,\rho}(t, mt^\beta \vec{v} )$ with $\varrho>1$ and $0<\beta<1$. 

To overcome partially such a difficulty, we introduce several auxiliary results to face this problem. In a first place, by the inverse Fourier transform formula of a radial function, we note that $u_{\alpha,\varrho}$ has the following representation
\begin{equation*}
u_{\alpha,\varrho}(t,x)=\frac{1}{\pi}\int_{0}^\infty E_\alpha(t^\alpha(1-\xi^{2\varrho}))\cos(x \xi)\, d\xi, \quad t>0, \quad x \in \RR,
\end{equation*}
where $C$ is a non-negative constant. 

The oscillatory behavior of $\xi \mapsto \cos(x \xi)$ in the preceding representation complicates the analysis of the asymptotic behavior of $u_{\alpha,\varrho}$. To address this, we state a generalization of \cite[Lemma 4.1]{Dipierro-Pellacci-Valdinoci-Verzini-2021}, which, in some way, separates the negative part of the solution from its positive part. We omit the proof since it is anologous to the case $\rho=1$ (see \cite[Lemma 4.1]{Dipierro-Pellacci-Valdinoci-Verzini-2021}).

\begin{lemma}\label{Lema:Leibniz}
 Let $\alpha\in(0,1)$ and $\varrho\ge1$. Then $u_{\alpha,\varrho}$ has the following representation 
    \begin{align*}
        u_{\alpha,\varrho}(t,x)=\frac{1}{\pi}\sum_{k=0}^{\infty}(-1)^ka_{k}(t,x),
    \end{align*}
where
\begin{align*}
    a_{0}(t,x)&=\int_{0}^{\frac{\pi}{2x}}E_{\alpha}(t^\alpha(1-\xi^{2\varrho}))\cos(x\xi)\, d\xi, \\
    a_k(t,x)&=(-1)^k\int_{\frac{\pi}{2x}(2k-1)}^{\frac{\pi}{2x}(2k+1)}E_{\alpha}(t^\alpha(1-\xi^{\varrho}))\cos(x\xi)\, d\xi, \quad k \geq 1.
\end{align*}
Further, the sequence $\{a_k(\cdot,\cdot)\}_{k\geq 1}$ satisfies the conditions of the \textit{Leibniz Criteria for alternating series}, that is, for every $(t,x) \in (0,\infty)\times \mathbb{R}$ and $k \geq 0$, the following statements hold:
\begin{enumerate}[$(a)$]
    \item $a_k(t,x)>0$,
    \item $\displaystyle \lim_{k \to \infty}a_k(t,x)=0$,  
    \item $a_{k+1}(t,x)<a_k(t,x)$, \quad for all $k \geq 1$, \quad \text{and}\quad $2a_0(t,x)>a_1(t,x)$.
\end{enumerate}
\end{lemma}

As a simple consequence of Lemma \ref{Lema:Leibniz} we obtain the following result (we left the proof to the reader).
\begin{lemma}\label{Lemma:Bounds:Leibniz}
   Let $\alpha\in(0,1)$ and $\varrho\ge1$. For all $m\in\NN$ we have 
    \begin{align*}
       u_{\alpha,\varrho}(t,x)&>\sum_{k=0}^{2m+1}(-1)^k a_{k}(t,x), \\
        u_{\alpha,\varrho}(t,x)&<\sum_{k=0}^{2m}(-1)^k a_{k}(t,x).
    \end{align*}
\end{lemma}
In particular,

\begin{align}\label{a0-a1}
    u_{\alpha,\varrho}(t,x)>a_{0}(t,x)-a_1(t,x).
\end{align}
%%%%%%%%%%%%%%%%%%%%% Lema Leibniz %%%%%%%%%%%%%%%%%%%%%%%%%%%%%%%%

Now, note that by means of the change of variable $r:=1-\xi^{2\varrho}$ in the expressions given in Lemma \ref{Lema:Leibniz}, the coefficients $a_0$ and $a_1$ can be rewritten as follows,
    \begin{align*}
        a_{0}(t,x)&=\frac{1}{2\varrho}\int_{1-(\frac{\pi}{2x})^{2\varrho}}^{1}E_{\alpha}(t^\alpha r)\cos(x\sqrt[2\varrho]{1-r}\,)(1-r)^{\frac{1-2\varrho}{2\varrho}}\,dr, \\
        a_{1}(t,x)&=-\frac{1}{2\varrho}\int_{1-(\frac{3\pi}{2x})^{2\varrho}}^{1-(\frac{\pi}{2x})^{2\rho}}E_{\alpha}(t^\alpha r)\cos(x\sqrt[2\varrho]{1-r}\,)(1-r)^{\frac{1-2\varrho}{2\varrho}}\,dr.
    \end{align*}

In what follows, we look for a lower estimate for $a_0$ and an upper estimate for $a_1$, in such a manner that $a_0(t,mt^\beta)-a_1(t,mt^\beta)\to\infty$ as $t\to \infty,$ for $\varrho\geq 1$ and $0<\beta<1$. Since $u_{\alpha,\varrho}$ satisfies \eqref{a0-a1}, this will imply that 
\[
\lim_{t\to\infty} u_{\alpha,\varrho}(t,mt^\beta)\to\infty,\quad \text{as}\quad t\to\infty,
\]
so the proof of Theorem \ref{Theorem:one:dimension} will be finished. For this purpose we prove the following auxiliary results.

\begin{lemma}\label{a0}
    Let $n\in\{2,3,\cdots\},$ $\rho\geq 1$ and $\alpha\in\big[\frac{1}{n},1)$. If $0<\ell<\frac{\pi}{2}<x$, then there is a function $c_0\colon \RR_+\times\RR_+\to \RR$ such that 
    \begin{equation*}
       a_0(t,x)\geq \frac{\cos(\ell)}{2\varrho}\left(\frac{x}{\ell}\right)^{2\varrho-1}\frac{\alpha}{t^{\alpha n}}\left[E_\alpha(t^\alpha)-E_\alpha\left(t^\alpha\left(1-\left(\frac{\ell}{x}\right)^{2\varrho}\right)\right)+c_0(t,x)\right], \quad t>0,
    \end{equation*}
    with
    \[
    \lim_{t \to \infty}\frac{c_0(t,mt^\beta)}{E_\alpha(t^\alpha)}=0,
    \]
    for all $\beta \in \mathbb{R}$ and $m>0$.
\end{lemma}

\bdemo
    Let $x>\frac{\pi}{2}$ and consider the function $h\colon [0,1)\to\RR$ defined by 
    \[
    h(r):= \cos(x\sqrt[2\rho]{1-r}\,)(1-r)^{\frac{1-2\varrho}{2\varrho}},\quad r\in [0,1).
    \]
    Since $\varrho\geq 1$, we have that $h$ is increasing in the interval $\Big(1-(\frac{\ell}{x})^{2\varrho},1\Big),$ with $\ell<\pi/2.$ Indeed, we note that
    \begin{align*}
        \frac{d}{dr}h(r)=\sin(x(1-r)^{\frac{1}{2\varrho}})\frac{x}{2\varrho}(1-r)^{\frac{1}{\varrho}-2}-\cos(x(1-r)^{\frac{1}{2\varrho}})\frac{1-2\varrho}{2\varrho}(1-r)^{\frac{1}{2\varrho}-2}.
    \end{align*}

Since $0<x(1-r)^{\frac{1}{2\varrho}}<\frac{\pi}{2}$, we have that,
    \begin{align*}
        0<\cos(x(1-r)^{\frac{1}{2\varrho}}) \,\,\, \text{and} \,\,\, 
        0<\sin(x(1-r)^{\frac{1}{2\varrho}}),
    \end{align*}
which implies that $h'(r)>0$ provided that $r\in\Big(1-(\frac{\ell}{x})^{2\varrho},1\Big)$. Also, by the choice of $\ell$, we have that $(1-(\frac{\ell}{x})^{2\varrho},1)\subset (1-(\frac{\pi}{2x})^{2\varrho},1)$. Therefore,
    \begin{align*}
        a_0(t,x)\geq& \frac{1}{2\varrho}\int_{1-(\frac{\ell}{x})^{2\varrho}}^1E_{\alpha}(t^\alpha r)\cos(x\sqrt[2\varrho]{1-r}\,)(1-r)^{\frac{1-2\varrho}{2\varrho}}\,dr \\
        &\geq \frac{\cos(\ell)}{2\varrho}\left(\frac{x}{\ell}\right)^{2\varrho-1}\int_{1-(\frac{\ell}{x})^{2\varrho}}^1 E_{\alpha}(t^\alpha r)dr.
    \end{align*}
It follows from Lemma \ref{Lemma:lower:estimate} that
\begin{align*}
    a_{0}(t,x)\geq& \frac{\cos(\ell)}{2\varrho}\left(\frac{x}{\ell}\right)^{2\varrho-1}\int_{1-(\frac{\ell}{x})^{2\varrho}}^1 \left(\frac{\alpha}{(t^\alpha r)^{n-1}}E_{\alpha}'(t^\alpha r)-\frac{1}{(t^\alpha r)^{n-1}}\sum_{k=0}^{[n-1+\frac{1}{2\alpha}]}\lambda_{k}(t^\alpha r)^k \right.\\ &\left.
    +\frac{1}{(t^\alpha r)^{n-1}}\sum_{k=n-1}^{[n-1+\frac{1}{2\alpha}]}\beta_{k}(t^\alpha r)^k\right)dr,\\
    \geq& \frac{\cos(\ell)}{2\varrho}\left(\frac{x}{\ell}\right)^{2\varrho-1}\int_{1-(\frac{\ell}{x})^{2\varrho}}^1\left(\frac{\alpha}{t^{\alpha(n-1)}}E_{\alpha}'(t^\alpha r)-\frac{1}{(t^\alpha r)^{n-1}}\sum_{k=0}^{[n-1+\frac{1}{2\alpha}]}\lambda_{k}(t^\alpha r)^k \right. \\ &\left. +\frac{1}{(t^\alpha r)^{n-1}}\sum_{k=n-1}^{[n-1+\frac{1}{2\alpha}]}\beta_{k}(t^\alpha r)^k\right)dr.
\end{align*}
%Procederemos con el cálculo de la primera integral. La aplicación del Teorema Fundamental del Cálculo en este paso es crucial para nuestro análisis, ya que nos proporcionará el término dominante $E_\alpha(t^\alpha)$. Este término será fundamental más adelante, ya que será responsable de la divergencia de $u_{\alpha,s}$.

By Fundamental Theorem of Calculus note that
\begin{align*}
    \int_{1-(\frac{\ell}{x})^{2\varrho}}^1\frac{\alpha}{t^{\alpha(n-1)}}E_{\alpha}'(t^\alpha r)\, dr=\frac{\alpha}{t^{\alpha n}}\left[E_\alpha(t^\alpha)-E_\alpha\left(t^\alpha\left(1-\left(\frac{\ell}{x}\right)^{2\varrho}\right)\right)\right].
\end{align*}
The other integral terms in the above lower estimate of $a_0(t,x)$ only produce terms of polynomial or logarithmic order, and they are dominated asymptotically by the function $E_\alpha(t^\alpha)$. Indeed, for $k\in\big\{0,1,\cdots, [n-1+\frac{1}{2\alpha}]\big\}$ we have that
\[
b_{k}(x):=
    \int_{1-(\frac{\ell}{x})^{2\varrho}}^1 r^{k+1-n} dr=\begin{cases} 
             \frac{1-(1-(\frac{\ell}{x})^{2\varrho}))^{k+2-n}}{k+2-n},& k+1-n \neq -1, \\
             \\ -\ln\big|1-(\frac{\ell}{x})^{2\varrho}\big|, & k+1-n = -1 .
             \end{cases}
\]             
%Observamos que la segunda rama tiene sentido considerarla, ya que $$M-2\leq M-1+\frac{1}{2\alpha} \implies M-2 \leq \left[M-1+\frac{1}{2\alpha}\right].$$
On the other hand, for $n-1 \leq k \leq [n-1+\frac{1}{2\alpha}]$ we have that 
\[
\int_{1-(\frac{\ell}{x})^{2\varrho}}^1 r^{k+1-n} dr=\frac{1-(1-(\frac{\ell}{x})^{2\varrho})^{k+2-n}}{k+2-n}.
\]
Thus, 
\begin{align*}
      %a_{0}(x,t)\geq& \frac{\cos(\ell)}{{2\varrho}}\left(\frac{x}{\ell}\right)^{2\varrho-1}\frac{\alpha}{t^{\alpha n}}\left[E_\alpha(t^\alpha)-E_\alpha\left(t^\alpha\left(1-\left(\frac{\ell}{x}\right)^{2\varrho}\right)\right)\right]\\
      %&-\frac{\cos(\ell)}{2\varrho}\left(\frac{x}{\ell}\right)^{2\varrho-1}\sum_{k=0}^{[n-1+\frac{1}{2\alpha}]}\lambda_{k}t^{\alpha k-\alpha(n-1)}b_{k}(t,x)\\
    %&+\frac{\cos(\ell)}{2\varrho}\left(\frac{x}{\ell}\right)^{2\varrho-1}\sum_{k=n-1}^{[n-1+\frac{1}{2\alpha}]}\gamma_{k}t^{\alpha k-\alpha(n-1)}\frac{1-(1-(\frac{\ell}{x})^{2\varrho})^{k+2-n}}{k+2-n}\\
 a_0(t,x)   \ge&\frac{\cos(\ell)}{2\varrho}\left(\frac{x}{\ell}\right)^{2\varrho-1}\frac{\alpha}{t^{\alpha n}}\left(E_\alpha(t^\alpha)-E_\alpha\left(t^\alpha\left(1-\left(\frac{\ell}{x}\right)^{2\varrho}\right)\right)-\frac{t^\alpha}{\alpha}\sum_{k=0}^{[n-1+\frac{1}{2\alpha}]}\lambda_{k}t^{\alpha k}b_{k}(x)\right.\\ &\left.
    +\frac{t^\alpha}{\alpha}\sum_{k=n-1}^{[n-1+\frac{1}{2\alpha}]}\beta_{k}t^{\alpha k}\frac{1-(1-(\frac{\ell}{x})^{2\varrho})^{k+2-n}}{k+2-n}\right) \\
    &=\frac{\cos(\ell)}{2\varrho}\left(\frac{x}{\ell}\right)^{2\varrho-1}\frac{\alpha}{t^{\alpha n}}\left[E_\alpha(t^\alpha)-E_\alpha\left(t^\alpha\left(1-\left(\frac{\ell}{x}\right)^{2\varrho}\right)\right)+c_0(t,x)\right],
\end{align*}
where
\begin{align*}
    c_0(t,x):=-\frac{t^\alpha}{\alpha}\sum_{k=0}^{[n-1+\frac{1}{2\alpha}]}\lambda_{k}t^{\alpha k}b_{k}(x)+
    \frac{t^\alpha}{\alpha}\sum_{k=n-1}^{[n-1+\frac{1}{2\alpha}]}\beta_{k}t^{\alpha k}\frac{1-(1-(\frac{\ell}{x})^s)^{k+2-{2\varrho}}}{k+2-n},
\end{align*}
and the constants $\beta_k$ and $\lambda_k$ have been defined in Lemma \ref{Lemma:lower:estimate}.
By the asymptotic behavior of Mittag-Leffler function (see Lemma \ref{asymp:mittag}), we have that
\begin{align*}
    \frac{c_0(t,mt^\beta)}{E_\alpha(t^\alpha)}= \frac{\displaystyle-\frac{t^\alpha}{\alpha}\sum_{k=0}^{[n-1+\frac{1}{2\alpha}]}\lambda_{k}t^{\alpha k}b_{k}(mt^\beta)+
    \frac{t^\alpha}{\alpha}\displaystyle\sum_{k=n-1}^{[n-1+\frac{1}{2\alpha}]}\beta_{k}t^{\alpha k}\frac{1-(1-(\frac{\ell}{mt^\beta})^s)^{k+2-n}}{k+2-n}}{\frac{1}{\alpha}e^t+o(1/t^{\alpha})}=o(1),\ \text{as}\ t\to\infty,
\end{align*}
which proves the result.
\fdemo

\begin{lemma}\label{a1}
    Let $n\in\{2,3,\cdots\},$ $\rho\geq 1$ and $\alpha\in\big[\frac{1}{n},1)$. If $\frac{3\pi}{2}<x$, then there is a function $c_1\colon \RR_+\times\RR_+\to \RR$ such that 
     \begin{align*}
         a_1(t,x)\leq \frac{\alpha}{2 t^\alpha\varrho }\left(\frac{2x}{\pi}\right)^{2\varrho-1} \left[E_\alpha\left(t^{\alpha}\left(1-\left(\frac{\pi}{2x}\right)^{2\varrho}\right)\right)-E_\alpha\left(t^{\alpha}\left(1-\left(\frac{3\pi}{2x}\right)^{2\varrho}\right)\right)+c_1(t,x)\right]
     \end{align*}
     with
     \begin{align*}
         \lim_{t \to \infty}\frac{t^{\alpha (n-1)} c_1(t,mt^\beta)}{E_\alpha(t^\alpha)}=0,
     \end{align*}
     for all $\beta \in \mathbb{R}$ and $m>0$.
\end{lemma}
\bdemo
    Since the function $r \mapsto (1-r)^{\frac{1-2\varrho}{2\varrho}}$ is increasing for $r\in (0,1)$ when $\varrho\geq 1$, we have that
    \begin{align*}
        a_1(t,x) %&\leq\int_{1-(\frac{3\pi}{2x})^s}^{1-(\frac{\pi}{2x})^s}E_{\alpha}(t^\alpha \rho)(1-\rho)^{\frac{1-s}{s}}\,d\rho \\
        &\leq \frac{1}{2\varrho}\left(\frac{2x}{\pi}\right)^{2\varrho-1}\int_{1-(\frac{3\pi}{2x})^{2\varrho}}^{1-(\frac{\pi}{2x})^{2\varrho}}E_{\alpha}(t^\alpha r)dr.
    \end{align*}
Further, since $x>\frac{3\pi}{2}$, we have that $1>\Big(\frac{3\pi}{2x}\Big)^{2\varrho}$. Thus, it follows from Lemma \ref{Lemma:upper:estimate} that
\begin{align*}
    a_1(t,x)&\leq \frac{1}{2\varrho}\left(\frac{2x}{\pi}\right)^{2\varrho-1}\int_{1-(\frac{3\pi}{2x})^{2\varrho}}^{1-(\frac{\pi}{2x})^{2\varrho}}(\alpha E'_{\alpha}(t^\alpha r)+\sum_{k=0}^{\left[\frac{3n-2}{2}\right]}\gamma_{k}(t^\alpha r)^k)\, dr \\
    &\leq  \frac{\alpha}{2t^\alpha \varrho}\left(\frac{2x}{\pi}\right)^{2\varrho-1} \left[E_\alpha\left(t^{\alpha}\left(1-\left(\frac{\pi}{2x}\right)^{2\varrho}\right)\right)-E_\alpha\left(t^{\alpha}\left(1-\left(\frac{3\pi}{2x}\right)^{2\varrho}\right)\right)\right]\\
    &+ \frac{1}{2\varrho}\left(\frac{2x}{\pi}\right)^{2\varrho-1}\sum_{k=0}^{[\frac{3n-2}{2}]}\gamma_{k}t^{\alpha k}\left(\frac{(1-(\frac{\pi}{2x})^{2\varrho})^{k+1}}{k+1}-\frac{(1-(\frac{3\pi}{2x})^{2\varrho})^{k+1}}{k+1}\right)\\
    &\leq  \frac{\alpha}{2t^\alpha \varrho}\left(\frac{2x}{\pi}\right)^{2\varrho-1} \left[E_\alpha\left(t^{\alpha}\left(1-\left(\frac{\pi}{2x}\right)^{2\varrho}\right)\right)-E_\alpha\left(t^{\alpha}\left(1-\left(\frac{3\pi}{2x}\right)^{2\varrho}\right)\right)+c_1(t,x)\right],
\end{align*}
where
\begin{align*}
    c_1(t,x):=\frac{t^\alpha}{\alpha}\sum_{k=0}^{[\frac{3n-2}{2}]}\gamma_{k}t^{\alpha k}\left(\frac{(1-(\frac{\pi}{2x})^{2\varrho})^{k+1}}{k+1}-\frac{(1-(\frac{3\pi}{2x})^{2\varrho})^{k+1}}{k+1}\right).
\end{align*}
By the asymptotic behavior of Mittag-Leffler function (see Lemma \ref{asymp:mittag}), we have
\begin{align*}
    \frac{t^{\alpha (n-1)} c_1(t,mt^\beta)}{E_\alpha(t^\alpha)}=\frac{\frac{t^{n\alpha}}{\alpha}\displaystyle\sum_{k=0}^{[\frac{3n-2}{2}]}\gamma_{k}t^{\alpha k}\left(\frac{(1-(\frac{\pi}{2mt^\beta})^{2\varrho})^{k+1}}{k+1}-\frac{(1-(\frac{3\pi}{2mt^\beta})^{2\varrho})^{k+1}}{k+1}\right)}{\frac{e^t}{\alpha}+o(1/t^{\alpha})}=o(1),\ \text{as}\ t\to\infty,
\end{align*}
which proves the result.
\fdemo 

\begin{theorem}\label{Teorema: diveregencia s>2}
    Let $d=1$, $n\in\{2,3,\cdots,\}$ and $\varrho \geq 1 $. For  all $\alpha \in [\frac{1}{n},1)$, $\beta \in (0,\frac{1}{2\varrho})$ and $m>0$, we have  
    \begin{align*}
        u_{\alpha,\varrho}(t,mt^\beta)>Cm^{2\varrho-1}t^{\beta(2\varrho-1)-\alpha n} E_{\alpha}(t^\alpha)[1-o(1)],
    \end{align*}
    where $C_1=\frac{\alpha\ell^{2-2\varrho}}{2\varrho}$ and $\ell \in (0,\frac{\pi}{2})$ is the unique solution of $\cos(\ell)=\ell$. In particular,
    \begin{align*}
        \lim_{t \to \infty}u_{\alpha,\varrho}(t,mt^\beta)=\infty.
    \end{align*}

\end{theorem}
\bdemo
Let \(\ell \in \big(0, \frac{\pi}{2}\big)\) be the unique solution of \(\cos(\ell) = \ell\), and \(x > \frac{\pi}{2}\). It follows from Lemma~\ref{Lemma:Bounds:Leibniz}, Lemma~\ref{a0}, and Lemma~\ref{a1} that
\begin{align*}
    u_{\alpha,\varrho}(t,x) \ge{}& \ \ell^{2 - 2\varrho} \frac{x^{2\varrho - 1} \alpha}{t^{\alpha n} 2\varrho} \Biggl[
    E_\alpha(t^\alpha) - E_\alpha\biggl(t^\alpha \bigl(1 - \bigl(\tfrac{\ell}{x}\bigr)^{2\varrho}\bigr)\biggr) + c_0(t,x)
    \Biggr] \\
    & + \frac{\alpha}{t^{\alpha} 2\varrho} \left(\frac{2x}{\pi}\right)^{2\varrho - 1} 
    \Biggl[
    E_\alpha\left(t^{\alpha}\left(1 - \left(\frac{3\pi}{2x}\right)^{2\varrho}\right)\right)
    - E_\alpha\left(t^{\alpha}\left(1 - \left(\frac{\pi}{2x}\right)^{2\varrho}\right)\right)
    - c_1(t,x)
    \Biggr] \\
    \ge{}& \ \frac{x^{2\varrho - 1}}{t^{\alpha n}} E_{\alpha}(t^\alpha) \Biggl[
    C_1 - C_1 \frac{E_\alpha\left(t^\alpha\left(1 - \left(\frac{\ell}{x}\right)^{2\varrho}\right)\right)}{E_{\alpha}(t^\alpha)} 
    + C_1 \frac{c_0(t,x)}{E_{\alpha}(t^\alpha)} \\
    & + C_2 \frac{t^{\alpha(n-1)} E_\alpha\left(t^\alpha \left(1 - \left(\frac{3\pi}{2x}\right)^{2\varrho}\right)\right)}{E_\alpha(t^\alpha)} 
    - C_2 \frac{t^{\alpha(n-1)} E_\alpha\left(t^\alpha \left(1 - \left(\frac{\pi}{2x}\right)^{2\varrho}\right)\right)}{E_\alpha(t^\alpha)} \\
    & - C_2 \frac{t^{\alpha(n-1)} c_1(t,x)}{E_\alpha(t^\alpha)}
    \Biggr],
\end{align*}
where 
\[
C_1 = \ell^{2 - 2\varrho} \frac{\alpha}{2\varrho} \quad \text{and} \quad C_2 = \frac{\alpha}{2\varrho} \left(\frac{2}{\pi}\right)^{2\varrho - 1}
\]
are positive constants.

Let us consider $x=mt^{\beta}$ for some large enough $t>0$. It follows from Lemma \ref{a0} and Lemma \ref{a1} that
\begin{align*}
    \frac{c_0(t,mt^\beta)}{ E_\alpha(t^\alpha)}=o(1),\quad \frac{t^{\alpha (n-1)}c_1(t,mt^\beta)}{E_{\alpha}(t^{\alpha})}=o(1),\quad t\to\infty.
\end{align*}
On the other hand, the asymptotic behavior of Mittag-Leffler function (Lemma \ref{asymp:mittag}) together with Bernoulli inequality imply that
\begin{align*}
    \frac{E_\alpha(t^\alpha\left(1-\ell^{2\varrho}m^{-2\varrho}t^{-2\varrho\beta})\right)}{ E_{\alpha}(t^\alpha)}&=\frac{e^{t((1-\ell^{2\varrho}m^{-2\varrho}t^{-2\varrho\beta})^\frac{1}{\alpha})}+o(1)}{e^t+o(1)}\\
    &\leq \frac{e^{t(1-\frac{1}{\alpha}\ell^{2\varrho}m^{-2\varrho}t^{-2\varrho\beta})}+o(1)}{e^t+o(1)}\\
    &=e^{-(\frac{1}{\alpha}\ell^{2\varrho}m^{-\varrho}t^{-2\varrho\beta+1})}+o(e^{-t}),
\end{align*}
as $t \to \infty$. Since Si $0<\beta <\frac{1}{2\varrho}$, we have that
\begin{align*}
     \frac{E_\alpha(t^\alpha\left(1-\ell^{2\varrho}m^{-2\varrho}t^{-2\varrho\beta})\right)}{ E_{\alpha}(t^\alpha)}=o(1),\quad t\to\infty.
\end{align*}
Analogously, for any $\varepsilon>0$ and $0<\beta<\frac{1}{2\varrho}$,
\begin{align*}
    \frac{t^{\alpha(n-1)}E_\alpha\left(t^\alpha\left(1-\varepsilon t^{-2\varrho\beta}\right)\right)}{ E_\alpha(t^{\alpha})}=t^{\alpha(n-1)}e^{-\varepsilon t^{-2\varrho\beta+1}}=o(1), \quad t\to\infty.
\end{align*}
Therefore, for all $0<\beta<\frac{1}{2\varrho}$ one gets 
\begin{align*}
    u_{\alpha,\varrho}(t,mt^\beta)\geq C_1 m^{2\varrho-1}t^{\beta(2\varrho-1)-\alpha n} E_{\alpha}(t^\alpha)[1-o(1)], \quad t\to\infty.
\end{align*}
In particular, by the asymptotic behavior of Mittag-Leffler function (Lemma \ref{asymp:mittag}),  $u_{\alpha,\varrho}(t,mt^{\beta})\to\infty$ as $t$ goes to infinity. 
\fdemo

\begin{remark}As we have mentioned before, in the particular case $n=1$ and $\varrho=1$, Theorem \ref{Teorema: diveregencia s>2} has been obtained by Dipierro, Pellacci, Valdinoci and Verzini \cite[Theorem 1.2]{Dipierro-Pellacci-Valdinoci-Verzini-2021}. 
\end{remark}

\section{Case $\varrho\in(0,1)$ revisited}\label{exponetial:invasion:speed}

As we have mentioned in the first section, the statement $i)$ in Theorem \ref{Theo:rho:general} suggests that a power-type function does not provide enough information for determining the optimal speed that distinguishes between the diverging and vanishing behavior of the fundamental solution of equation \eqref{Main:Equation} in the case $\varrho \in (0,1)$. To overcome this difficulty, we consider the function $\theta_{m,\beta}(t) = (e^{m t^{\beta}} - 1)$ for some $m > 0$ and $\beta > 0$. The subtraction of $1$ allows to the function $\theta$ to start from $0$ at $t=0$. By choosing appropriate values for $m$ and $\beta$, we can potentially control the speed and manner in which $\theta_{m,\beta}$ grows, allowing for a more refined understanding of the behavior of solutions to the equation in question. This approach may lead to insights that cannot be easily obtained by using power-type functions alone, particularly when dealing with equations involving fractional powers ($\varrho \in (0,1)$ in this case).

%\begin{theorem}\label{Theorem:rho:in:(0,2):again}
%Let $d\in\NN$, $\alpha\in(0,1)$, and $\vec{v}\in\mathbb{S}^{d-1}$. Let $u_{\alpha,\varrho}$ be the solution of \eqref{Main:Equation:3} with $\varrho\in(0,1)$. Define $\theta_{m,\beta}(t) = (e^{m t^{\beta}} - 1)$ for some $m > 0$ and $\beta > 0$. The following statements hold.
%\begin{enumerate}[$(i)$]
%\item If $\beta\in(0,1)$ and $m>0$ then 
%\[
%\lim_{t\to\infty}u_{\alpha,\varrho}(t,\theta_{m,\beta}(t)\vec{v})=\infty.
%\]

%\item If $\beta>1$ and $m>0$ then %=\begin{cases}\infty, 
%\[
%\lim_{t\to\infty}u_{\alpha,\varrho}(t,\theta_{m,\beta}(t)\vec{v})=0.
%\]
%\item If $\beta=1$, and $0<m\le \frac{1-\gamma_\alpha}{d+2\varrho}$, where $\gamma_\alpha=(1-\alpha)\alpha^{\frac{\alpha}{1-\alpha}}$, then 
%\[
%\lim_{t\to\infty}u_{\alpha,\varrho}(t,\theta_{m,\beta}(t)\vec{v})=\infty.
%\]
%\item If $\beta=1$ and $m>\frac{1}{d+2\varrho}$, then 
%\[
%\lim_{t\to\infty}u_{\alpha,\varrho}(t,\theta_{m,\beta}(t)\vec{v})=0.
%\]
%\end{enumerate}
%\end{theorem}

\subsection*{Proof of Theorem \ref{Theorem:rho:in:(0,2)}} According to Lemma \ref{Lemma:Subordination}, the fundamental solution of \eqref{Main:Equation} is given by 
\[
u_{\alpha,\varrho}(t,x)=t^{-\alpha}\int_0^\infty e^{s}\ \mu^\varrho_s(x)W_{-\alpha,1-\alpha}(-t^{-\alpha}s)ds, \quad t> 0, \quad x\in\RR^d,
\]
where $(\mu^\varrho_s)_{s\ge 0}$ is the convolution semigroup associated to the symbol $\xi\mapsto \|\xi\|^\varrho$. It is well known, see, e.g. \cite[Section 3]{Jacob-2001}, that the Fourier transform of $\mu^\varrho_s$ is given by 
\[
\widetilde{\mu}^{\varrho}_s(\xi)=e^{-s\|\xi\|^{2\varrho}},\quad s>0, \ \xi\in\RR^d.
\]
As we have mentioned before, there is not closed formulas for $\mu_s^{\varrho}$ except in the cases $\varrho=1$ and $\varrho=\frac{1}{2}$. However, in any case there are positive constants $C_1$ and $C_2$ such that 
\begin{equation}\label{Lower:Upper:Poison}
C_1\frac{t}{\Big(\|x\|^2+t^{\frac{1}{\varrho}}\Big)^{\frac{d+2\varrho}{2}}}\le\mu_{t}^\varrho(x)\le C_2  \frac{t}{\Big(\|x\|^2+t^{\frac{1}{\varrho}}\Big)^{\frac{d+2\varrho}{2}}}, \qquad x \in \R^n, \quad t>0.
\end{equation}

By a similar argument to the proof of Theorem \ref{Theorem:rho=2} \((i)\), the positivity and the lower bounds of \( \mu_t^\varrho \) given in \eqref{Lower:Upper:Poison} imply that
\begin{align*}
 u_{\alpha,\varrho}(t,x) \gtrsim \frac{te^{t}}{((2t)^{1/\varrho}+\|x\|^2)^{\frac{d+2\varrho}{2}}} W_{-\alpha,1}(-t^{1-\alpha})(1-o(1)), \qquad t\to \infty.
\end{align*}
Considering \( x = \theta_{m,\beta}(t)\vec{v} \) and the asymptotic behavior of the Wright functions, we have that 
\begin{align*}
    u_{\alpha,\varrho}(t,x)  \gtrsim \frac{t^{\frac{1}{2}}e^{t-\gamma_\alpha t}}{((2t)^{1/\varrho}+(e^{mt^\beta}-1)^2)^{\frac{d+2\varrho}{2}}}(1-o(1)), \qquad t \to \infty.
\end{align*}
Moreover, taking \( t > 0 \) sufficiently large such that \( e^{mt^\beta} > e^{mt^\beta} - 1 > (2t)^{1/\varrho} \), the last inequality implies that
\begin{align*}
     u_{\alpha,\varrho}(t,x)  \gtrsim t^{\frac{1}{2}} e^{t(1 - t^{\beta-1}m(d+2\varrho) - \gamma_\alpha)}(1-o(1)), \qquad t \to \infty.
\end{align*}
Therefore, it follows that if \( \beta < 1 \), then 
\begin{align*}
    \lim_{t \to \infty} u_{\alpha,\varrho}(t,\theta_{m,\beta}(t)\vec{v}) = \infty. 
\end{align*}

Furthermore, in the case \( \beta = 1 \), we have  $\lim_{t \to \infty} u_{\alpha,\varrho}(t,\theta_{m,\beta}(t)\vec{v}) = \infty$, provided that 
\[
1 - m(d + 2\varrho) - \gamma_\alpha \geq 0,
\]
which is equivalent to the condition $0<m\le \frac{1-\gamma_\alpha}{d+2\varrho},$ and the statement $(iii)$ follows.

For parts \((ii)\) and \((iv)\), we consider the upper bound of the inequalities \eqref{Lower:Upper:Poison} and the Taylor expansion of the exponential function, to get
\begin{align*}
    u_{\alpha,\varrho}(t,x) &\lesssim t^{-\alpha} \int_{0}^{\infty} e^s \frac{s}{(s^{1/\varrho} + \|x\|^2)^{\frac{d+2\varrho}{2}}} W_{-\alpha,1-\alpha}(-t^{-\alpha}s) \, ds \\
    &\leq \frac{t^{-\alpha}}{\|x\|^{d+2\varrho}} \int_{0}^{\infty} \sum_{k=0}^{\infty} \frac{s^{k+1}}{k!} W_{-\alpha,1-\alpha}(-t^{-\alpha}s) \, ds.
\end{align*}
We denote \( \Phi(r) := W_{-\alpha,1-\alpha}(-r) \). Using Lemma \ref{asymp:mittag} and Proposition \ref{Moments}, we have that
\begin{align*}
        u_{\alpha,\varrho}(t,x) &\lesssim \frac{t^{-\alpha}}{\|x\|^{n+2\varrho}} \int_{0}^\infty \sum_{k=0}^\infty \frac{s^{k+1}}{k!} \Phi(st^{-\alpha}) \, ds = \frac{t^{-\alpha}}{\|x\|^{d+2\varrho}} \sum_{k=0}^\infty \frac{1}{k!} \int_{0}^\infty s^{k+1} \Phi(st^{-\alpha}) \, ds \\
                &= \frac{1}{\|x\|^{d+2\varrho}} \sum_{k=0}^\infty \frac{(t^\alpha)^{k+1}}{k!} \int_{0}^\infty u^{k+1} \Phi(u) \, du = \frac{1}{\|x\|^{d+2\varrho}} \sum_{k=0}^\infty \frac{(t^\alpha)^{k+1}}{k!} \frac{\Gamma(k+2)}{\Gamma(\alpha(k+1)+1)} \\
                &= \frac{1}{\|x\|^{d+2\varrho}} \sum_{k=0}^\infty \frac{(t^\alpha)^{k+1}}{k!} \frac{(k+1)!}{\alpha(k+1) \Gamma(\alpha(k+1))} = \frac{1}{\alpha \|x\|^{d+2\varrho}} \sum_{k=0}^\infty \frac{(t^\alpha)^{k+1}}{\Gamma(\alpha k + \alpha)}\\
                &= \frac{t^{\alpha}}{\alpha \|x\|^{d+2\varrho}} E_{\alpha,\alpha}(t^\alpha) = \frac{t^{\alpha}}{\alpha \|x\|^{d+2\varrho}} \left( \frac{1}{\alpha} t^{1-\alpha} e^{t} + O(t^{2\alpha}) \right).
\end{align*}
If \( x = (e^{m t^{\beta}} - 1 )\vec{v}\), for \( \beta > 1 \) and some unitary vector $\vec{v}$, it follows that
\begin{align*}
        \lim_{t \to \infty} u_{\alpha,\varrho}(t,(e^{m t^\beta}-1)\vec{v}) = 0.
\end{align*}
For the case $\beta =1,$ the statement $(iv) $ follows providing that $m(d+2\varrho)>1$.

\begin{remark} Note that, similar to Theorem \ref{Theorem:rho=2}, in the case $\beta=1$, there is a gap where the asymptotic behavior of $u_{\alpha,\varrho}(t, (e^{mt}-1)\vec{v})$ remains unknown. Unlike the previous scenario, this gap depends not only on the parameter $\alpha$ but also on the dimension $d$ and on the order $\varrho$, where the gap decreases as $d$ increases, as well as, $\alpha \to 1^{-}$. In such a case, we have that 

\[
\lim_{t \to \infty} u_{\alpha,\varrho}(t, (e^{mt}-1) \vec{v}) = 
\begin{cases}
\infty, & m \leq \frac{1-\gamma_\alpha}{d+2\varrho}, \\
0, & m > \frac{1}{d+2\varrho}.
\end{cases}
\]

While our method does not provide information about the asymptotic behavior of \( u_{\alpha,\varrho}(t,(e^{mt}-1)\vec{v}) \) for \( m \) within the range 

\[
\frac{1-\gamma_\alpha}{d+2\varrho} < m \leq \frac{1}{d+2\varrho},
\] 
this gap closes as \( \alpha \to 1^{-} \). 
Thus, in the limite case, we recover the classical result

\[
\lim_{t \to \infty} u_{1,\varrho}(t,(e^{mt}-1)\vec{v}) =
\begin{cases}
\infty, & m \leq \frac{1}{d+2\varrho}, \\
0, & m > \frac{1}{d+2\varrho}.
\end{cases}
\] 

\fdemo

\end{remark}

\begin{comment}
    \begin{remark}Note that, similar to Theorem \ref{Theorem:rho=2}, in the case $\beta=1$, there is a gap where the asymptotic behavior of $u_{\alpha,\varrho}(t, (e^{mt}-1)\vec{v})$ remains unknown. Unlike the previous scenario, this gap depends not only on the parameter $\alpha$ but also on the dimension $d$ and on the order $\varrho$, diminishing as $d$ and $\rho$ increases.

\end{remark}

\end{comment}

%From this point, the proof of $(i)$, $(ii)$ and $(iii)$ follows the same ideas of the proof of Theorem \ref{Theorem:rho=2}. \textcolor{red}{*} For $(iv)$, we observe that

\bibliographystyle{acm}
\bibliography{Biblio.bib}

\end{document}